\documentclass{article}
\usepackage{amsfonts, amsmath, amsthm, amssymb}
\usepackage[all]{xy}

\newtheorem{theorem}{Theorem}[section]
\newtheorem{lemma}[theorem]{Lemma}
\newtheorem{cor}[theorem]{Corollary}
\newtheorem{prop}[theorem]{Proposition}

\newtheorem{definitiontemp}[theorem]{Definition}
\newenvironment{definition}{\begin{definitiontemp}
\normalfont}{\end{definitiontemp}}

\theoremstyle{remark}
\newtheorem*{remark}{Remark}

\newcommand{\cl}{\mathcal{C}\ell^{\bullet}}
\newcommand{\akg}{A_\kappa (\mathfrak{g})}
\newcommand{\exten}{\text{Ext}}
\newcommand{\pair}{(\cdot\,,\cdot)}
\newcommand{\gk}{\hat{\mathfrak{g}}_\kappa}
\newcommand{\gr}{\mathcal{G}r}
\newcommand{\gt}{G(\mathcal{O})}
\newcommand{\frakk}{\mathfrak{k}}
\newcommand{\fraki}{\mathfrak{i}}
\newcommand{\gl}{\Check{G}}
\newcommand{\hl}{\Check{H}}
\newcommand{\nk}{\mathfrak{n}(\mathcal{K})}
\newcommand{\hk}{\hat{\mathfrak{h}}_\kappa}
\newcommand{\fl}{\mathcal{F}\ell}
\newcommand{\lk}{\mathcal{L}_{\kappa}}

\newcommand{\vac}{\left|0\right>}
\newcommand{\height}{\text{ht}}
\newcommand{\Hk}{S_{\kappa}(G)}

\newcommand{\vkg}{V_\kappa (\mathfrak{g})}
\newcommand{\Nk}{N(\mathcal{K})}

\newcommand{\g}{\mathfrak{g}}
\newcommand{\V}{\mathcal{V}}

\newcommand{\hcrit}{\hat{\mathfrak{h}}_{\kappa-\kappa_c}}
\newcommand{\C}{\mathbb{C}}
\newcommand{\Waff}{W_{\text{aff}}}
\newcommand{\bfrak}{\mathfrak{b}}
\newcommand{\h}{\mathfrak{h}}
\newcommand{\hla}{\left|\lambda\right\rangle}
\newcommand{\hchi}{\left|\chi\right\rangle}

\newcommand{\bla}{b^{\bar{\lambda}}(w)}
\newcommand{\mv}{Mirkovi\'{c}-Vilonen }
\newcommand{\coh}{H^{\infty/2+\bullet}}
\newcommand{\cohr}{H^{\bullet}}
\newcommand{\cohdr}{H^{\bullet}_{DR}}
\newcommand{\lb}{\mathcal{L}}

\newcommand{\iw}{\fraki}
\newcommand{\dett}{\text{det}}
\newcommand{\dimm}{\text{dim}}

\newcommand{\Iw}{I}

\newcommand{\lam}{\lambda_w}
\newcommand{\wba}{\bar{w}}
\newcommand{\ho}{\mathfrak{h}(\mathcal{O})}
\newcommand{\dind}[2]
{\genfrac{}{}{0pt}{}{#1}{#2}}
\newcommand{\K}{\mathcal{K}}
\newcommand{\Oo}{\mathcal{O}}
\newcommand{\n}{\mathfrak{n}}
\newcommand{\crit}{\kappa_{c}}
\newcommand{\m}{\mathfrak{m}}

\newcommand{\A}{\mathcal{A}}
\newcommand{\skg}{S_\kappa(G)}
\newcommand{\akh}{A_\kappa(\h)}
\newcommand{\mods}{\text{mod}}
\newcommand{\gO}{\g(\Oo)}
\newcommand{\nO}{\n(\Oo)}
\newcommand{\nkp}{\nk_{\dagger}}
\newcommand{\nop}{\nO_{\dagger}}
\newcommand{\Nkp}{\Nk_{\dagger}}
\newcommand{\determinant}{\text{det}}
\newcommand{\dimension}{\text{dim}}
\newcommand{\gm}{\mathbb{G}_m}
\newcommand{\N}{\mathcal{N}}

\begin{document}

\title{The BRST reduction of the chiral Hecke algebra}
\author{Ilya Shapiro%
\thanks%
{The main results in this paper are part of the author's doctoral
dissertation, written at the University of Chicago under the
supervision of Alexander Beilinson. The author wishes to thank
Professor Beilinson for many illuminating discussions over the years
as well as careful readings of the drafts of this paper. Additional work on this text was carried out at UC Davis, MPI (Bonn), IHES, and the University of Waterloo.  Many thanks are also due to Roman Bezrukavnikov for helpful discussions.} }

\date{}
\maketitle

\begin{abstract}

We explore the relationship between de Rham and Lie algebra
cohomologies in the finite dimensional and affine settings.  In
particular, given a $\gk$-module that arises as the global sections
of a twisted $D$-module on the affine flag manifold, we show how to
compute its untwisted BRST reduction modulo $\nk$ using the de Rham
cohomology of the restrictions to $\Nk$ orbits.  A similar
relationship holds between the regular cohomology and the Iwahori
orbits on the affine flag manifold. As an application of the above,
we describe the BRST reduction of the chiral Hecke algebra as a
vertex super algebra.

\end{abstract}


\section{Introduction.}
A way of looking at geometric representation theory is as an attempt
to match up algebraic objects that naturally arise in the study of
representations of groups or algebras, with geometric objects which
are perhaps easier to study.  An early example of this is the
Borel-Weil-Bott theorem that constructs irreducible representations
of a reductive group via the sheaf cohomology of equivariant line
bundles on the flag manifold of the group.

Expanding on the above approach, one may obtain representations of a
Lie algebra $\g$ by considering the global sections of a $D$-module
on a homogeneous space of $G$, where $\g=Lie(G)$.  In the case of a
reductive group $G$ and its flag manifold $G/B$, we obtain in this
way all representations of $\g$ with the trivial central character.
This is part of the work of Beilinson-Bernstein \cite{bb81} which
was aimed at proving the Kazhdan-Lusztig conjecture.\footnote{The
Kazhdan-Lusztig conjecture was independently demonstrated around the
same time by Brylinski-Kashiwara in \cite{brylkash2} using very
similar methods.}  Note that $G$ itself acts on both categories, via the twisting of $\g$-modules by the adjoint action of $g\in G$ and the pullback of $D$-modules along the action of $g\in G$ on $G/B$.  This identification is compatible with the actions, thus we can say that
these two categories are but two incarnations of the correct
analogue of the representations of $G$ on the space of ``functions"
on $G/B$. For more on this point of view see \cite{newbook}.

A more complete version of the result of \cite{bb81} is that
representations of $\g$ with other central characters may be
obtained from appropriately twisted $D$-modules on $G/B$ with the
twisting corresponding to the central character. Thus the center  of
$U\g$, i.e., the center of the enveloping algebra of $\g$, serves as
the space of ``spectral parameters" for a decomposition of its
category of representations. It coincides with the Bernstein center
of the category of $\g$-representations. This type of ``spectral
decomposition" of the category of $\g$-representations and the
identification of the ``fibers" with categories of geometric origin
has proven itself to be very useful.

Pursuing this further, we can also consider the setting of affine
Kac-Moody Lie algebras $\gk$, which are infinite dimensional
analogues of reductive Lie algebras.  Here $\g$ is as above and
$\gk$ is a central extension\footnote{The central extensions are parameterized by the levels $\kappa$ which are invariant inner products on the Cartan subalgebra (see Sec. \ref{notation} for more details).  In this paper, outside of the introduction, we are only concerned with the negative integral level.} of the Lie algebra of the loop group
$G(\K)$ best thought of as parameterizing maps from the punctured
formal disc $D^\times$ to $G$.  The story becomes more interesting
at this point and links up with the geometric Langlands program.
Namely the space of ``spectral parameters" now called local
Langlands parameters is the moduli stack parameterizing de Rham $\gl$-local
systems on $D^\times$ (i.e., $\gl$ principal bundles on $D^\times$ with an automatically flat connection), where we denote by $\gl$ the Langlands dual
group of $G$.\footnote{A good reference for the notion of a category over a stack is \cite{stack}.} Thus to each local Langlands parameter $\chi$ one must
attach an appropriate subcategory $\gk-\mods_\chi$ of $\gk$-modules
that is stable under the action of $G(\K)$. Considerable progress
has been made in this direction by Frenkel-Gaitsgory in the case of
a critical level $\kappa$; it is well surveyed in \cite{newbook}.
The key aspect of this particular value of the level is that the
center is very large. There is a conjecture of Beilinson, stated in
the introduction to \cite{torus} that addresses the case of the
negative integral level. The chiral Hecke algebra $\skg$ of
Beilinson-Drinfeld plays a central role there.

In this paper we restrict our attention to the geometric part of the picture above. It
fits into the general framework as follows: one considers only the
subcategory of $\gk$-modules with ``support" in the substack of
regular singular connections on $D^\times$ with nilpotent residue; these  have a conjectural interpretation as $D$-modules.\footnote{This is known as the tamely ramified case of the local geometric Langlands conjecture.  Furthermore, we will focus almost exclusively on the unramified case.}  This substack can be described concretely as the quotient stack $\N/\gl$, where $\N$ is the nilpotent cone of $\gl$.  Briefly, the elements of $\N$ represent the residue of the connection and the quotient by $\gl$ accounts for the gauge transformations.  In fact the category of $D$-modules on $\fl$ is itself naturally a category over the stack $\N/\gl$ (this fact underlies \cite{upgrade, upgrade2}).  This explains the appearance of the regular singular condition on the connection.

Basically we want to emulate the correspondence between
$D_{G/B}$-modules and $\g$-modules with the trivial central
character.\footnote{We point out at this time that what is actually
accomplished here is only a first, albeit important, step.}  The
role of $D$-modules is still played by $D$-modules, now on the
affine flag manifold, however there is no obvious candidate for the
subcategory of $\g$-modules specified by the triviality of the
central character, as the enveloping algebra of $\gk$, or rather its
appropriate analogue, has no center. This corresponds to the fact
that the moduli stack of de Rham local systems, discussed above,  has no non-constant
global functions \cite{torus}. The notion of ``support" replaces the center, and the meaning to the ``support" of $\gk$-modules is given through the consideration of the categories of modules over the twists of the chiral Hecke algebra by $\gl$-local systems. In the case under consideration, i.e. the analogue of the trivial central character, we look at local systems with regular singularities and nilpotent residue.\footnote{This is a general principle of generating representations of the smaller vertex algebra $\vkg$ by considering twisted representations of the larger $\skg$ that contains it.  By keeping track of the twisting we obtain a measure of control over the representations of $\vkg$ that we allow.  Here and below when speaking about $\gk$-modules we are implicitly using the fact that they are canonically identified with $\vkg$-modules, where $\vkg$ is the Kac-Moody vertex algebra.}

In short, one wants to interpret $D$-modules on the affine flags as certain special $\gk$-modules just as in the Beilinson-Bernstein localization
theorem. Currently there
are partial analogues of the localization
theorem in the context of the negative integral level. Namely, in the work of Beilinson-Drinfeld \cite{bd} and
Frenkel-Gaitsgory \cite{localization} it is shown that the modules
arising from appropriately twisted $D$-modules on either the affine
flag manifold $\fl$ or the affine Grassmannian $\gr$ embed into the
category of $\gk$-modules for the $\kappa$ corresponding to the
twisting. However, identifying the image of the above embedding is
problematic. As mentioned above, the main candidate for the space of spectral
parameters, namely the center of the enveloping algebra, that was
used for this purpose in the finite dimensional case, is absent
here. Instead one should, conjecturally, use the chiral Hecke
algebra $\Hk$. We postpone any discussion of $\skg$ to Sec.
\ref{chabrief} and need only point out that $\skg$ is obtained from
the twisted global sections of a $D$-module on the affine
Grassmannian (denoted by $\widetilde{\Oo}_{\gl}$), it is
$\gl$-equivariant (so that it can be twisted by a $\gl$-local system), and $\skg^{\gl}$ is the Kac-Moody vertex algebra
$\vkg$ whose representation theory is the same as that of $\gk$.

The localization conjecture for the affine
flag manifold that serves as one of the motivations for the present
paper is an important special case of the conjecture outlined in the
introduction of \cite{torus} (where this more general conjecture is
settled for the considerably simpler commutative case). Namely, it
is conjectured that there is an equivalence between, roughly
speaking, the category of appropriately twisted equivariant representations of
$\Hk$ and the product of several copies of the category of
$D$-modules on $\fl$. More precisely, consider the following
commutative diagram of functors
$$
\xymatrix{(\skg_{\N},\gl)-\mods\ar[rd]|{\Gamma(\N,-)^{\gl}} &  &
D_{\fl}-\mods\ar[ld]|{\Gamma(\fl,-\otimes\lb_{\kappa+\chi})}
\ar[ll]_-{F_\chi}\\
 & \gk-\mods &
}
$$  and a few words of explanation for the symbols used are in order.  As is repeatedly mentioned above, $\skg$ is a $\gl$-equivariant vertex algebra so that any $\gl$-local system $\phi$ on $D^\times$ gives rise to a new chiral algebra on the punctured disc that we denote by $\skg_\phi$.  Recall that the moduli stack of $\gl$-local systems with a regular singularity and nilpotent residue is given by $\N/\gl$ so that we obtain in this way $\skg_{\N}$, a bundle of chiral algebras on $\N$ that is $\gl$-equivariant.   We denote by $(\skg_{\N},\gl)-\mods$ the category of $\gl$-equivariant
$\skg_{\N}$-modules, with the notation for the other two categories being self explanatory.  The functor $F_\chi$ is based on the concepts of \cite{functorz} and \cite{upgrade}.  It is roughly $\Gamma(\fl,(\mathcal{Z}(\Oo_{\gl})\star-)\otimes\lb_{\kappa+\chi})$ where $\mathcal{Z}$ is the functor from \cite{functorz}\footnote{This functor maps $Rep\, \gl$, the category of representations of $\gl$, to $D^I_{\fl}-\mods$, the category of Iwahori equivariant $D$-modules on $\fl$.  Here $\Oo_{\gl}$ is the $\gl$-module of functions on $\gl$.} and $\star$ is the fusion product, see \cite{fusion, functorz}
for example.  We say roughly because it is not clear how it lands in the $\skg_{\N}$-modules, to see this one needs some ideas of \cite{upgrade}.  Let us mention that $G(\K)$ acts on each category in the diagram and the functors commute with this action.  There is an action of $Rep\,\gl$ on both sides of $F_\chi$, obvious on the left and via $\mathcal{Z}$ on the right, and $F_\chi$ commutes with it.

The top arrow becomes (conjecturally) an equivalence
of categories if we sum over appropriate representatives $\chi$.  Namely we consider the affine Weyl group dot action on the weight
lattice of $G$, with respect to the level $\kappa-\crit$, and $\chi$ is the only dominant regular, in the affine sense, weight in a given orbit.\footnote{This is the complementary point of view to our notion of sufficiently negative level, discussed in a Remark following Lemma \ref{deltasemi}.}
The conjecture  solves the problem of
identifying precisely which $\gk$-modules come from $D$-modules on
the affine flags: they are the ones that extend to an equivariant action of $\skg_{\N}$, the $\gl$-equivariant
bundle of chiral algebras on $\N$ that contains $\vkg$.

The above is the tamely ramified case of the conjecture in \cite{torus}.  Let us now consider the unramified case.  It concerns the category $D_{\fl}^{m.a.}-\mods$ of monodromy annihilators, i.e. $D$-modules $M$ on $\fl$ such that the monodromy action of \cite{functorz} on $\mathcal{Z}(V)$, with $V$ any representation of $\gl$, becomes trivial on $\mathcal{Z}(V)\star M$.  When restricted to this subcategory the functor $F_\chi$ lands in $(\skg,\gl)-\mods$ which is the subcategory of $(\skg_{\N},\gl)-\mods$ supported at $0\in\N$.  This is because the lack of monodromy ensures that the action vertex operators of $\skg$ are no longer multi-valued and so we do not need to twist it by local systems in order to get rid of this complication. Thus we obtain the following diagram:

$$
\xymatrix{(\skg,\gl)-\mods\ar[rd]|{-^{\gl}} &  &
D^{m.a.}_{\fl}-\mods\ar[ld]|{\Gamma(\fl,-\otimes\lb_{\kappa+\chi})}
\ar[ll]_-{F_\chi}\\
 & \gk-\mods &
}
$$ and the conjecture is that $F_\chi$ is an equivalence after summing over $\chi$ as above.

It would be interesting to investigate the relationship between $D_{\gr}-\mods$ and $D^{m.a.}_{\fl}-\mods$. Namely as seen in the following diagram:
$$
\xymatrix{(\skg_{\N},\gl)-\mods  & & &
D_{\fl}-\mods\ar[lll]_{\quad\quad F_{2\rho}}
\\
(\skg,\gl)-\mods\ar[u]^{i_\bullet} & & & D_{\gr}-\mods\ar[u]^{\pi^*}\ar[lll]_{\quad\quad\Gamma(\gr,(\widetilde{\Oo}_{\gl}\star-)\otimes\lb_\kappa)}
}
$$ the $D$-modules from $\gr$ provide a large supply of monodromy annihilators via $\pi^*$.  However it is not difficult to come up with, using \cite{upgrade}, examples of m.a. $D$-modules that are not pulled back from $\gr$, at least not via $\pi^*$.\footnote{For $G=PGL_2$ there is an automorphism $\sigma$ of $\fl$ such that $\sigma^*\pi^*$ also produces monodromy annihilators.}  One may naively conjecture, based on a similar result \cite{flfromgr} on the level of derived categories, that the category $D^{m.a.}_{\fl}-\mods$ is obtained from $D_{\gr}-\mods$ via base change from the stack $0/\gl$ to $0/\Check{B}$.  This would have an interesting consequence that a $\gl$-equivariant $\skg$-module can be twisted not only by a $\gl$-representation, but more generally, by a $\Check{B}$-representation.

Another question is how to describe the image of $F_{2\rho}\circ\pi^*$ above.  A possible answer is discussed after Corollary \ref{conj} and involves the BRST functor.  The key is that the BRST reduction of an $\skg$-module that comes directly from a $D$-module on the affine Grassmannian has a very compact and conjecturally characterizing form in terms of the de Rham cohomology of the restrictions of the original $D$-module to the semi-infinite orbits in $\gr$.

In this paper we compute by a mixture of algebraic and geometric
methods, $\coh(\nk,\skg)$, i.e., the semi-infinite cohomology of
$\nk$ with coefficients in $\skg$.  It follows from general
considerations that as $\skg$ is a vertex algebra, so is
$\coh(\nk,\skg)$ and we explicitly describe its vertex algebra
structure.

For a $\gk$-module $M$, the motivation for considering its BRST
reduction, as $\coh(\nk, M)$ is called, lies in the suggestive fact
that it is a $\hcrit$-module.  Thus the problem shifts from the
domain of the non-commutative $G$ to the more accessible case of its
commutative torus $H$.\footnote{In fact in \cite{torus} the conjecture is checked  for the case of $G=H$.}  Broadly described, a
possible approach to the problem of identifying $D^{m.a.}_{\fl}-\mods$ with
$(\skg,\gl)-\mods$ consists of first trying to enumerate the images
of the objects on both sides under appropriate functors and then
hope to lift this identification to the original categories.  The
functors are, to first approximation,
$\oplus_{w\in\Waff}\cohdr(S_w,i_w^! -)$ on the $D$-module
side\footnote{The $S_w$ are the $\Nk$-orbits which are labeled by
the elements of the affine Weyl group $\Waff$.} and the BRST
reduction on the $\skg$ side.\footnote{This is illustrated by Corollaries \ref{conj} and \ref{conj2}.}  The latter requires a few words of
explanation.  If $M$ is a $\gl$-equivariant $\skg$-module then by
the results of this paper $\coh(\nk, M)$ is a $\gl$-equivariant
module over the vertex algebra of global sections of a
$\gl$-equivariant bundle of vertex algebras over $\gl/\hl$.  The
category of such modules is then equivalent to the category of
$\hl$-equivariant modules over the vertex algebra that is the fiber
of the bundle above over $\hl\in\gl/\hl$.  This fiber is an
enlargement of the much studied lattice Heisenberg vertex
algebra,\footnote{Coincidentally, the lattice Heisenberg vertex
algebra is the chiral Hecke algebra for $G=H$.} in fact the lattice
Heisenberg vertex algebra is exactly its cohomological degree $0$
part.\footnote{The remaining part is roughly $H^\bullet(\n,\C)$.}
The representation theory of the lattice Heisenberg is well
understood; it has a finite number of irreducible modules. This
observation agrees with the fact that one should really consider a
product of several copies of $D_{\fl}-\mods$ as corresponding to
$(\skg,\gl)-\mods$.  Recalling now that we still have the
$\hl$-grading and the remaining part of the fiber vertex algebra, we
see that one has roughly the same type of object as
$\oplus_{w\in\Waff}\cohdr(S_w,i_w^! M)$.  We hope that the results
and methods of this paper will provide a way to illuminate the
relationship between the algebraic, i.e., the representation
theoretic side and the geometric, i.e., the $D$-module side of the
above conjectural correspondence.

This text is organized as follows. Section \ref{lieandderham}
contains comparison theorems between Lie algebra and De Rham
cohomologies that we will subsequently need. Some of the results in
this section (in particular the ones pertaining to the finite
dimensional situation) are believed to be part of the folklore;
unfortunately we cannot cite a reference other than this text. The
proofs provided here are based on A. Voronov's semi-infinite
induction (alternatives are demonstrated in the finite case).  It is
worth noting that Theorems \ref{fincoh} and \ref{fincoh2} are
essentially the same, but the proofs illustrate very different
approaches. They address the finite-dimensional case. Theorem
\ref{semicoh} is the main theorem of this section, it deals with the
semi-infinite version of the infinite dimensional case.

Section \ref{brstreduction} is devoted to the computation of the
BRST reduction of the chiral Hecke algebra, first as a module over
the Heisenberg Lie algebra (Corollary \ref{hmodofreduction}), and
finally, in the main theorem of the paper (Theorem
\ref{maintheorem}), as a vertex algebra.  Some of the ingredients
used in the proof are Theorem \ref{semicoh} and the \mv theorem
\cite{satake,satake2}.

In section \ref{appendix} we provide some auxiliary information that
the reader should find useful. Namely, a brief overview of the
Beilinson-Drinfeld chiral Hecke algebra is included (Sec.
\ref{chabrief}). No references containing a construction were available for
citation, however a brief discussion can be found in \cite{thebook}. The language of the highest weight algebras is introduced
(Sec. \ref{hwaappendix}) as it is useful for stating the main
results of the paper.  Also included in Sec. \ref{hwaappendix} are
certain details on how the Heisenberg Lie algebra module
structure on a vertex algebra determines the vertex operators modulo
the knowledge of the highest weight algebra.

Some sources containing the background material for this paper that
we recommend are \cite{cha, string, thebook}. Finally, the terms
semi-infinite cohomology and BRST reduction are used interchangeably
and we refer the reader to \cite{voronov2} for the definitions. A
sketch of the relevant details is given in the discussion preceding
Lemma \ref{deltasemi}. The matter of notation
is addressed below. Since our sources do not use
mutually compatible notation, we made some choices that are to the
best of our knowledge consistent.

\subsection{Some notational conventions.}\label{notation}

We encourage the reader to quickly skim this section and to review whenever necessary.

The group $G$ that we consider is a simple algebraic group over
$\C$, and $\g$ is its Lie algebra. Some of the groups and
algebras that we need are  the Lie algebras $\bfrak\subset\g$, the Borel
subalgebra, $\n=[\bfrak,\bfrak]$, the nilpotent subalgebra, and $\h=\bfrak/\n$ the Cartan Lie algebra; the corresponding groups are denoted by $B$, $N$, and $H$.  We reserve $\bfrak^-$, $\n^-$, etc for the opposite versions, i.e., $\n^-$ is the sum of the negative root spaces. Note that $\h$ is sometimes used to
denote a subalgebra of $\bfrak$ but this requires a choice, the same
holds for $H$.

Put $\K=\C((t))$ and $\Oo=\C[[t]]$. Let $\g(\K)=\g\hat{\otimes}\K$,
and define $\g(\Oo)$, $\n(\K)$, and $\n(\Oo)$ similarly. Denote by
$G(\K)$ and $N(\K)$ the algebraic loop groups of $G$ and $N$, by
$G(\Oo)$ and $N(\Oo)$ the subgroups of positive loops. Denote by
$\nkp$ and $\nop$, $\n(\K)\oplus t\h(\Oo)$ and $\n(\Oo)\oplus
t\h(\Oo)$ respectively.

Given an invariant inner product $\pair_{\kappa}$ on $\g$, the
affine Kac-Moody Lie algebra $\gk$ is defined as the central
extension $\g(\K)^{\sim}$ of $\g(\K)$, with the cocycle $\phi$ given
by $\phi(x\otimes f,y\otimes g)=-(x,y)_{\kappa} \text{Res} f dg$,
where $x, y\in\g$ and $f, g\in \mathcal{K}$. For the purposes of
this paper $\pair_\kappa=\kappa\pair_0$ with $\kappa < -h^\vee$,
where $\pair_0$ is the normalized invariant inner product on $\g$
(i.e., $(\theta,\theta)_0=2$, where $\theta$ is the highest weight
of the adjoint representation) and $h^\vee=1+(\rho,\theta)_0$
($\rho=\frac{1}{2}\sum_{\alpha>0}\alpha$) the dual Coxeter number of
$\g$. This ensures that the level super line bundle is defined and
twisting by it makes the global sections functor exact and faithful
\cite{localization}. We note that
$\pair_{\crit}=-\frac{1}{2}\pair_{Kil}=-h^\vee\pair_0$.

Let $\Gamma$ denote the co-weight lattice and $\Check{\Gamma}$ the
weight lattice of $G$.  Write $\lb_\chi$ for the line bundle with
total space $G\times_{B}\C_{-\chi}$ for $\chi\in\Check{\Gamma}$. We point out that for us $\C_\chi$ denotes a non-trivialized line on which $\h$ (or $\Check{\h}$) acts via the character (or co-character) $\chi$. As usual $W$
and $\Waff$ denote the Weyl and the affine Weyl groups respectively,
note that $\Waff=\Gamma\rtimes W$.  In the finite setting the dot
action of $W$ is defined by $w\cdot\chi=w(\chi+\rho)-\rho$, where
$\rho$ is the half sum of the positive roots.  In the affine setting
the dot action depends on the level $\kappa$ and for $w\in \Waff$
with $w=\lambda_w \bar{w}$, is given by
$w\cdot\chi=\bar{w}\cdot\chi-\kappa(\lambda_{w})$.

To emphasize the role of $\rho$ we use the convention that $\chi$ is
called dominant if $(\chi + \rho)(H_\alpha) \notin \{-1,-2,-3,...\}$
for each positive coroot $H_\alpha$.  We say that $\chi$ is dominant
regular if $\chi-\rho$ is dominant.  We note that it is very common
to call the latter dominant, we do not follow that convention.

Denote by $I$ the Iwahori subgroup of $G(\Oo)$, more precisely,
$I=\text{ev}^{-1}(B)$ where $\text{ev}: G(\Oo)\rightarrow G$ is the
usual evaluation map. Set $I^+=\text{ev}^{-1}(N)$.  We will use
$\fraki$, $\fraki^+$ for $Lie(I)$, $Lie(I^+)$ respectively.  Let
$\fl$ denote the affine flag manifold and $\gr$ the affine
Grassmannian, roughly speaking $\fl=G(\K)/I$ and $\gr=G(\K)/G(\Oo)$.

We reserve $\mathcal{T}_X$ and $\Oo_X$ for the sheaves of vector
fields and functions on $X$ respectively.  If $C^\bullet$ is a
complex, then $C^\bullet[n]$ denotes a degree shift, i.e., the
degree $k$ component of $C^\bullet[n]$ is $C^{k+n}$.

\section{Lie algebra and De Rham cohomologies.}\label{lieandderham}
We are interested in reducing the Lie algebra cohomology (usual or
semi-infinite) computations for modules that arise geometrically
as twisted global sections of a $D$-module on a certain $G$-space,
to the computation of the De Rham cohomology of the $D$-module
itself restricted to orbits.  We begin with the motivational
finite-dimensional setting and proceed to the case of interest,
the affine setting.

\subsection{The finite-dimensional setting.}\label{finite}
Let $X$ be an homogeneous $G$ space.  Then by differentiating the
$G$ action we obtain a map of Lie algebras
$\alpha:\g\longrightarrow\Gamma(X,\mathcal{T}_X)$, which after
taking the dual gives
$\Gamma(X,\Omega_X^i)\longrightarrow\bigwedge^i
\g^*\otimes\Gamma(X,\Oo_X)$.  Furthermore if $M$ is a
left\footnote{All $D$-modules in the finite setting are left by
default, though we consider the right ones in Theorem
\ref{fincoh2}.  In the affine setting, only the right $D$-modules
exist.} $D$-module on $X$, then $\Gamma(X,M)$ is a $\g$ module,
and we have a map
$\Gamma(X,M\otimes\Omega_X^{\bullet})\longrightarrow\bigwedge^{\bullet}
\g^*\otimes\Gamma(X,M)$.  If $X$ is affine, the complex on the
left computes $\cohdr(X,M)$, while the one on the right computes
$\cohr(\g,\Gamma(X,M))$, and as our map commutes with the
differentials, it descends to the cohomology, namely
$$\alpha^*:\cohdr(X,M)\longrightarrow\cohr(\g,\Gamma(X,M)).$$

In addition, if the action of $G$ on $X$ extends to that of $G'$
in which $G$ is normal then both sides above are $\g'/\g$-modules
and the map is compatible with this action.  Note that
$\cohdr(X,M)$ is a trivial $\g'/\g$-module.  Furthermore, in the
case when $X$ is a $G$ torsor $\alpha^*$ is an isomorphism even on
the level of complexes.

Let us apply this observation to the following situation.  Given a
$D$-module $M$ on $G/B$, one may consider
$\cohr(\mathfrak{n},\Gamma(G/B, M \otimes \lb_{\chi}))$ as a
$\h$-module. We should immediately restrict our attention to $\chi$
dominant regular\footnote{See the remark following Theorem
\ref{fincoh} for the non dominant regular $\chi$ case.} as this
ensures the exactness of $\Gamma(G/B,-\otimes\lb_{\chi}))$. In that
case we have the following:

\begin{theorem}\label{fincoh}
Let $M$ be a D-module on $G/B$, and $X_w \subset G/B$ the $N$ orbit
labeled by $w \in W$, let $\chi$ be dominant regular, then as
$\h$-modules
$$\cohr(\mathfrak{n},\Gamma(G/B, M \otimes
\lb_{\chi}))\cong\bigoplus_{w \in W}\cohdr (X_{w}, i^{!}_{w}M)
\otimes \C_{w \cdot (-2\rho - \chi)}.$$
\end{theorem}

\begin{remark}
The proof given below, while illuminating, is ultimately a
digression.  The reader may skip to Theorem \ref{fincoh2} which,
along with its proof, is a baby version of the one in the affine
setting.
\end{remark}

\begin{proof}
Recall that we have a notion of length for the elements $w$ of the Weyl group $W$, in particular the length $\ell(w)$ is equal to the dimension of the corresponding $N$-orbit $X_w$. We observe that $G/B$ has a filtration (see the next paragraph) $S_i=\coprod_{\ell(w)\leq
i} X_w$ which equips $M$ with a filtration in the derived category
with associated graded factors $i_{w*}i^!_w M$.  Applying
$\Gamma(G/B,-\otimes\lb_{\chi}))$, we get a filtration on
$\Gamma(G/B, M \otimes \lb_{\chi})$.  This reduces the theorem to
the special case of $M=i_{w*}M_0$ as the $\h$ action on the
cohomology of the factors is then different for different $w$'s and
so the spectral sequence degenerates and
$\cohr(\mathfrak{n},\Gamma(G/B, M \otimes \lb_{\chi}))$ is
canonically isomorphic to $\bigoplus_{w \in
W}\cohr(\mathfrak{n},\Gamma(G/B, i_{w*}i^!_w M \otimes
\lb_{\chi}))$.

Since the above type of argument is used repeatedly in the rest of the paper we provide some additional details.  Consider a decomposition of a space $Y=\coprod Z_i$ with $T_n=\coprod_{i\leq n}Z_i$ closed in $Y$.  Let $\iota_n:T_n\hookrightarrow Y$, $\alpha_n:T_{n-1}\hookrightarrow T_n$, $j_n:Z_n \hookrightarrow T_n$ and $i_n: Z_n\hookrightarrow Y$ so that $i_n=\iota_n\circ j_n$ and $\iota_{n-1}=\iota_n\circ\alpha_n$.  If $M$ is a $D$-module on $Y$ then by considering $\iota_n^! M$ on $T_n$ and the decomposition $T_n=T_{n-1}\coprod Z_n$ we obtain a distinguished triangle in the derived category $\alpha_{n*}\alpha_n^!\iota_n^! M\rightarrow \iota_n^! M\rightarrow j_{n*}j_n^!\iota_n^! M.$  If we apply $\iota_{n*}$ to it we get another distinguished triangle $$\iota_{(n-1)*}\iota_{n-1}^! M\rightarrow \iota_{n*}\iota_{n}^! M\rightarrow i_{n*}i^!_{n} M.$$  The filtration on $M$ is thus given by $\iota_{n*}\iota^!_n M$ with the associated graded factors $i_{n*}i^!_n M$.  Since $\Gamma(G/B,-\otimes\lb_\chi)$ is exact by \cite{bb81}, it preserves distinguished triangles and applying the cohomological functor $H^\bullet(\n,-)$ we obtain the desired spectral sequence.  In fact one can avoid any reference to the machinery of spectral sequences by using induction and long exact sequences that will be canonically split exact using the $\h$ action.

We are ready to proceed, begin with $w_0$, the longest element in $W$, i.e., the element corresponding to the big cell in $G/B$.  Dropping the
subscript in $M_0$, we have
\begin{align*}
\cohr(\mathfrak{n},\Gamma(G/B, i_{w_0*}M \otimes
\lb_{\chi}))&\cong\cohr(\mathfrak{n},\Gamma(X_w, M \otimes
\lb_{\chi}|_{X_{w_0}}))\\
&\cong\cohr(\mathfrak{n},\Gamma(X_w, M ))\otimes
\C_{w_0(-\chi)}\\
&\cong\cohdr(X_{w_0},M)\otimes \C_{w_0(-\chi)},
\end{align*}
the last step follows from the discussion above as $X_{w_0}$ is an
$N$ torsor.  Note that we do not need $\chi$ to be dominant
regular here.

To prove the theorem for other $w\in W$ we reduce to the case of
$w_0$ using the following observation.  Let $Y_w$ in $G/B\times G/B$
be the $G$ orbit through $(B,wB)$, denote by $p_1$ and $p_2$ the
restriction to $Y_w$ of the projections onto the factors. For $M$ a
$D$-module on $G/B$, set $\widetilde{M}^w=p_{2*}p_1^*M$, then
$$\Gamma(G/B,M\otimes
\lb_{\chi})\cong R\Gamma(G/B,\widetilde{M}^w\otimes
\lb_{w^{-1}\cdot\chi})$$ as $\g$-modules.  Let us suppress the
exponent in $\widetilde{M}^w$ once it is established which $w$ we
are using.

Now let $M$ be a $D$-module on $X_w$. Consider the diagram

$$
\xymatrix{G/B & Y_{w^{-1}w_0}\ar[l]_{p_1}\ar[r]^{p_2} & G/B\\
X_w\ar[u]^{i_w} & Y_{w^{-1}w_0}'\ar[l]_{p_1'}\ar[r]^{p_2'}\ar[u]^i
& X_{w_0}\ar[u]^{i_{w_0}}}
$$  where the left square above is Cartesian by definition (i.e., $Y'$ is defined by the diagram itself), $p_1'$ has affine space fibers, and
$p_2'$ is an isomorphism.  So that $p_{2*}p_1^*i_{w*}M\cong
p_{2*}i_*p_1'^*M\cong i_{w_0 *}p_{2*}'p_1'^*M$, hence
$i_{w_0}^!\widetilde{i_{w*}M}^{w^{-1}w_0}\cong p_{2*}'p_1'^*M$. The
proof is then completed by the following chain of isomorphisms:
\begin{align*}
\cohr(\mathfrak{n},\Gamma(G/B, i_{w*}M \otimes \lb_{\chi}))&\cong
\cohr(\mathfrak{n},R\Gamma(G/B, \widetilde{i_{w*}M} \otimes \lb_{w_0
w\cdot\chi}))\\&\cong\cohr(\mathfrak{n},\Gamma(X_{w_0},
i^!_{w_0}\widetilde{i_{w*}M}))\otimes
\C_{w\cdot(-2\rho-\chi)}\\
&\cong\cohdr(X_{w_0},i_{w_0}^!\widetilde{i_{w*}M})\otimes
\C_{w\cdot(-2\rho-\chi)}\\
&\cong\cohdr(X_{w_0},p_{2*}'p_1'^*M)\otimes
\C_{w\cdot(-2\rho-\chi)}\\
&\cong\cohdr(Y_{w^{-1}w_0}',p_1'^*M)\otimes \C_{w\cdot(-2\rho-\chi)}\\
&\cong\cohdr(X_w,M)\otimes \C_{w\cdot(-2\rho-\chi)}.
\end{align*}
\end{proof}
\begin{remark}
The assumption that $\chi$ be dominant regular is necessary,
however there is a way to replace $M$ by $\widetilde{M}^w$, very
similar to the method used in the proof of Theorem \ref{fincoh} in
such a way that we have for $w^{-1}\cdot\chi$ dominant regular
$$R\Gamma(G/B,M\otimes\lb_{\chi})
\cong\Gamma(G/B,\widetilde{M}^w\otimes\lb_{w^{-1}\cdot\chi}).$$ When
$\chi$ isn't dominant regular but $w^{-1}\cdot\chi$ is\footnote{Such
a $w\in W$ exists if and only if
$\left<\Check{\alpha},\chi+\rho\right>\neq 0$ for all $\alpha\in
R$.}, this reduces the problem to our familiar case. The
construction of $\widetilde{M}^w$ is immediate from the observation
that for any character $\chi$, we have that
$\Gamma(G/B,i_{e*}\Oo_e\otimes\lb_{\chi})$ is the Verma module with
highest weight $-2\rho-\chi$, while for $\chi$ dominant regular
$\Gamma(G/B,i_{w!}\Oo_w\otimes\lb_{\chi})$ is the Verma module with
highest weight $w\cdot(-2\rho-\chi)$. Explicitly we set
$\widetilde{M}^w=p_{2!}p_1^* M$, where $p_1$ and $p_2$ are as in the
proof of Theorem \ref{fincoh}.  This ``intertwining functors"
construction originates in \cite{bb83}.
\end{remark}

\emph{So far we have been  using left $D$-modules implicitly, however in
the affine setting only right $D$-modules exist, and so we switch to
using them exclusively at this point.} Furthermore the proof of
Theorem \ref{fincoh} does not immediately generalize to that
setting. It was included because its very geometric nature appealed
to us. Now we must switch to a more algebraic approach that directly
generalizes. We begin with some preliminaries.

The following is a version of the Shapiro Lemma\footnote{It has a
semi-infinite analogue \cite{voronov}.}. Note that $\g$ must be finite dimensional.
For a finite dimensional $V$, we use $\dett(V)$ to denote its top
exterior power $\wedge^{\dimm(V)}V$; it is a non-trivialized line.

\begin{lemma}\label{shapiro}
Let $\frakk\subset\g$ be a Lie subalgebra, then there is a natural
isomorphism: $$H^\bullet(\frakk,
M\otimes\determinant(\g/\frakk)^*)\stackrel{\sim}{\longrightarrow}
H^\bullet(\g,Ind_\frakk^\g M)[dim\g-dim\frakk]$$ where $M$ is a
$\frakk$-module.
\end{lemma}
\begin{proof}
If $L$ is a finite dimensional Lie algebra and $N$ an $L$-module,
then there is an isomorphism
$$H^\bullet(L,N\otimes\determinant(L))\stackrel{\sim}{\longrightarrow}
H_\bullet(L,N)[-\dimension(L)]$$ given by the contraction of
$\determinant(L)$ with forms $\omega\in\wedge^\bullet L^*$.  One
checks that the map commutes with the differentials and it is
clearly an isomorphism on the level of complexes.  There is a map in
the other direction obtained by
$$H_\bullet(L,N)=H_\bullet(L,(N\otimes\determinant(L))\otimes\determinant(L^*))\rightarrow
H^\bullet(L,N\otimes\determinant(L))[\dimension(L)]$$ similarly
through contraction. The following chain of isomorphisms completes
the proof:

\begin{align*}
H^\bullet(\frakk,
M\otimes\determinant(\g/\frakk)^*)&\stackrel{\sim}{\rightarrow}
H_\bullet(\frakk,
M\otimes\determinant(\g/\frakk)^*\otimes\determinant(\frakk^*))[-\dimension(\frakk)]\\
&\cong H_\bullet(\frakk,
M\otimes\determinant(\g^*))[-\dimension(\frakk)]\\
\intertext{by Shapiro Lemma} &\stackrel{\sim}{\rightarrow}
H_\bullet(\g, \text{Ind}_\frakk^\g
(M\otimes\determinant(\g^*)))[-\dimension(\frakk)]\\
\intertext{By universality there is a map of $\g$-modules
$\text{Ind}_\frakk^\g (M\otimes\determinant(\g^*))\rightarrow
\text{Ind}_\frakk^\g (M)\otimes\determinant(\g^*)$ that is
compatible with the natural filtration on the modules and is an
isomorphism on the associated graded pieces.  Thus it is an
isomorphism of modules:} &\stackrel{\sim}{\rightarrow} H_\bullet(\g,
\text{Ind}_\frakk^\g
(M)\otimes\determinant(\g^*))[-\dimension(\frakk)]\\
&\stackrel{\sim}{\rightarrow} H^\bullet(\g,\text{Ind}_\frakk^\g
M)[\dimension\g-\dimension\frakk].
\end{align*}

Note that the map of the Lemma can be written down explicitly as
follows.  Observe that $\determinant(\g/\frakk)^*$ is naturally a
line in $\wedge^\bullet\g^*$ and whereas there is no canonical map
from $\wedge^\bullet\frakk^*$ to $\wedge^\bullet\g^*$, there is one
from $\wedge^\bullet\frakk^*\otimes\determinant(\g/\frakk)^*$ to
$\wedge^\bullet\g^*[\dimension\g-\dimension\frakk]$.  Tensoring this
map with $M\hookrightarrow \text{Ind}_\frakk^\g M$ yields the
required isomorphism.
\end{proof}

By a (right) $D$-module of delta functions at $x$ in $X$ we mean the $D$-module $i_{x*}\C$ where $i_x: \{x\}\hookrightarrow X$; we denote it by $\delta_x$.  If $x=B\in G/B$ then the sections of $\delta_x$ (as a $\g$-module) can be described explicitly as $U\g/U\bfrak=\text{Ind}_{\bfrak}^{\g}\C$.  Recall that such an object is called a Verma module.  In the representation theory of $\g$ one also has a Co-Verma module $\text{Coind}_{\bfrak^-}^{\g} \C$ and everything in between called semi-induced modules \cite{voronov, voronov2}.  The precise definition of the semi-induced module is not straightforward, however what we need is the fact that all of these $\g$-modules have the same character, i.e., agree as $\h$-modules.  At least as $\n$-modules they can be constructed through co-induction followed by induction (see the proof of the Lemma below).  Furthermore, each is well adapted to a particular (co)homology theory.  More precisely, $H_\bullet(\n^-,\text{Ind}_{\bfrak}^{\g}\C)=\C$, $H^{\bullet}(\n,\text{Coind}_{\bfrak^-}^{\g} \C)=\C$, etc.

When we consider appropriately twisted delta functions, i.e., $\delta_x\otimes\lb_\chi$ for $\chi$ sufficiently dominant, then (see \cite{kk}) the $\g$-module $\Gamma(G/B, \delta_x\otimes\lb_\chi)$ is simple and thus its cohomology may be computed by identifying it with any one of the (isomorphic in this case) semi-induced modules. This is the idea behind the proof of the Lemma below as well as its affine analogues.

\begin{lemma}\label{delta} Let $\delta_x$ be the right $D$-module of delta
functions at $x\in G/B$ and $\chi-2\rho$ be dominant regular, then
$$H^\bullet(\n, \Gamma(G/B, \delta_x\otimes\lb_\chi))\cong\lb_\chi|_x\otimes
\determinant(\n/s_\n x)^*[-\dimension(\n/s_\n x)]$$ where $s_\n x$
is the stabilizer in $\n$ of $x\in G/B$.
\end{lemma}

\begin{proof}
We note that it is sufficient to prove this statement for $x=wB$ for
$w\in W$, because for every $y\in G/B$, $\Gamma(G/B,
\delta_y\otimes\lb_\chi)$ is a twist of one of $\Gamma(G/B,
\delta_{wB}\otimes\lb_\chi)$ by an element of $N$.

Observe that $\Gamma(G/B, \delta_{wB}\otimes\lb_\chi)$ is a simple
$\g$-module. So we can identify it with a semi-induced module of
Voronov \cite{voronov}.  As a result\footnote{Alternatively, we can
obtain the same result by transferring the $D$-module from $X_w$ to
$X_{w_0}$ as in the proof of Theorem \ref{fincoh}.} we obtain a
description of $\Gamma(G/B, \delta_{wB}\otimes\lb_\chi)$ as
$\text{Ind}_{\n\cap \n_w}^\n\text{Coind}_0^{\n\cap \n_w}
\lb_\chi|_{wB}$ as an $\n$-module, where $\n_w=w\n w^{-1}$.  The
Lemma is then a consequence of the following isomorphisms that are
versions of Shapiro Lemma:
\begin{align*}
H^\bullet(0,\lb_\chi|_{wB}\otimes\determinant(\n/\n\cap
\n_w)^*)&\stackrel{\sim}{\leftarrow}
H^\bullet(\n\cap \n_w,\text{Coind}_0^{\n\cap \n_w}\lb_\chi|_{wB}\otimes\determinant(\n/\n\cap \n_w)^*)\\
&\stackrel{\sim}{\rightarrow}H^\bullet(\n,\text{Ind}_{\n\cap
\n_w}^\n\text{Coind}_0^{\n\cap \n_w}
\lb_\chi|_{wB})[\dimension(\n/s_\n x)],
\end{align*}
where the second isomorphism is Lemma \ref{shapiro}. Note  the use
of triviality of the $\determinant(\n/\n\cap \n_w)^*$ as an
$\n\cap\n_w$-module.

\end{proof}

The Corollary below is a consequence of the identification of
homology and cohomology explained in the proof of Lemma
\ref{shapiro}.
\begin{cor}\label{deltahom} With the assumptions of Lemma \ref{delta},
$$H_\bullet(\n, \Gamma(G/B,
\delta_x\otimes\lb_\chi))\cong\lb_\chi|_x\otimes \determinant(s_\n
x)[\dimension(s_\n x)].$$
\end{cor}

\begin{remark}
We follow the convention that dictates that the Lie algebra homology
is placed in negative degrees.  More precisely, $H_{-i}(\n, M)$ is a subquotient of $\wedge^i\n\otimes M$.
\end{remark}

\begin{theorem}\label{fincoh2}
Let $M$ be a right D-module on $G/B$, and $\chi-2\rho$ be dominant
regular, $n=\text{dim}(\n)$, then as $\h$-modules
$$H_\bullet(\mathfrak{n},\Gamma(G/B, M \otimes
\lb_{\chi}))\cong\bigoplus_{w \in W}\cohdr (X_{w}, i^{!}_{w}M)
\otimes \C_{-w \cdot (\chi-2\rho)}[n-\ell(w)].$$
\end{theorem}
\begin{remark}
We follow the convention that dictates that the De Rham cohomology
is placed in both positive and negative degrees. More precisely, the left exact functor $\Gamma$ is applied to the complex $\wedge^{-i}\mathcal{T}_X\otimes M$ that is confined to the negative degrees.
\end{remark}

\begin{proof}
As in the proof of Theorem \ref{fincoh} we may reduce to a
$D$-module of the form $i_{w*}M$ for some $w\in W$ and $M$ a
$D$-module on $X_w$.  The action of $N$ on $X_w$ yields the
following short exact sequence:
$$\mathcal{S}tab_w\stackrel{\alpha}{\longrightarrow} \Oo_{X_w}\otimes\n\stackrel{\beta}{\longrightarrow} \mathcal{T}_{X_w}$$
where $\mathcal{S}tab_w$ is the kernel of the action map $\beta$.
Choose a section $s$ of $\beta$, define
$$\psi:\bigwedge\nolimits^i\mathcal{T}_{X_w}\otimes\text{det}
(\mathcal{S}tab_w)\rightarrow\bigwedge\nolimits^{i+n-\ell(w)}\n\otimes\Oo_{X_w}$$
by $\omega\otimes v\mapsto s(\omega)v$, note that $\psi$ does not
depend on the choice of $s$. Then $\psi$ extends to
$$\widetilde{\psi}:i_{w\cdot}(M\otimes\lb_{\chi}|_{X_w}\otimes
\bigwedge\nolimits^\bullet\mathcal{T}_{X_w}\otimes\text{det}(\mathcal{S}tab_w))[n-\ell(w)]\rightarrow
i_{w*}M\otimes\lb_{\chi}\otimes\bigwedge\nolimits^{\bullet}\n$$
where $\widetilde{\psi}$ is a morphism of complexes of sheaves on
$G/B$ that we intend to show is actually a quasi-isomorphism (after
passing to $R\Gamma$ it yields the isomorphism of the Theorem).

Since the $N$-action trivializes both $\text{det}(\mathcal{S}tab_w)$
and $\lb_{\chi}|_{X_w}$ they contribute only a twist by a
$\h$-character and we have a map on the cohomologies: $$\cohdr
(X_{w}, M) \otimes \C_{-w \cdot (\chi-2\rho)}[n-\ell(w)]\rightarrow
H_\bullet(\mathfrak{n},\Gamma(G/B, i_{w*}M \otimes \lb_{\chi})).$$

Since $M$ has a finite resolution by finite sums of
$\mathcal{D}_{X_w}$ and their direct summands, we may assume that
$M=\mathcal{D}_{X_w}$.  In this case both sides are finite
dimensional vector bundles over $X_w$ and the map is a morphism of
$\Oo_{X_w}$-modules.  Over $x\in X_w$,  the map becomes $$\cohdr
(X_{w}, \delta_x) \otimes \C_{-w \cdot
(\chi-2\rho)}[n-\ell(w)]\rightarrow
H_\bullet(\mathfrak{n},\Gamma(G/B, i_{w*}\delta_x \otimes
\lb_{\chi})),$$ which is an isomorphism by Corollary \ref{deltahom}.
This completes the proof.
\end{proof}

The statement of Theorem \ref{fincoh2} is made in terms of Lie
algebra homology because the proof to us seemed most natural in that
case (it avoids relative determinants for now), however it can be
easily reformulated in terms of cohomology, namely
$$\cohr(\mathfrak{n},\Gamma(G/B, M \otimes
\lb_{\chi}))\cong\bigoplus_{w \in W}\cohdr (X_{w}, i^{!}_{w}M)
\otimes \C_{w \cdot (- \chi)}[-\ell(w)],$$ compare this with the remark
following Theorem \ref{semicoh}.

\begin{remark}
It was pointed out to us by A. Beilinson that Theorem \ref{fincoh}
(and thus Theorem \ref{fincoh2}) can be obtained as a consequence of
the Beilinson-Bernstein localization theorem.  The disadvantage of
course is that while the proof given below is very simple, it uses
both the center of $U\g$, and the localization theorem of
Beilinson-Bernstein; neither is present in the affine case.

\begin{proof} Observe that
\begin{align*}
\cohr(\mathfrak{n},\Gamma(G/B, M \otimes \lb_{\chi})) &\cong
\exten^{\bullet}_{U\n}(\C,\Gamma(G/B, M \otimes
\lb_{\chi}))\\
&\cong \exten^{\bullet}_{U\g}(U\g\otimes_{U\n}\C,\Gamma(G/B, M
\otimes \lb_{\chi})).
\end{align*}

Since $\Gamma(G/B, M \otimes \lb_{\chi})$ is a module obtained from
a twisted $D$-module, the center $\mathcal{Z}(U\g)$ of $U\g$ acts on
it via $\phi(-w_0(\chi))$, where $\phi:\h^*\twoheadrightarrow
\text{Spec} \mathcal{Z}(U\g)$ is the Harish-Chandra map.
Furthermore it acts in the same way on the Verma modules $\{V(w\cdot
(-w_0(\chi)))|w\in W\}=\{V(w\cdot(-2\rho-\chi))|w\in W\}$. Note that
$U\g\otimes_{U\n}\C$ on the other hand is a superposition of all
Verma modules, which as a sheaf on $\text{Spec} \mathcal{Z}(U\g)$ is
locally free near $\phi(-w_0(\chi))$ as $\phi$ is \'{e}tale there.

Let $\m_{\chi}$ be the maximal ideal in $\mathcal{Z}(U\g)$
corresponding to $\phi(-w_0(\chi))$. Let $F_{\bullet}$ be the
forgetful functor from the category of $U\g/\m_{\chi}$-modules to
the category of $U\g$-modules.  It admits an obvious left adjoint
$F^*$, namely the restriction to $\phi(-w_0(\chi))\in
\text{Spec}\mathcal{Z}(U\g)$.  The following chain of isomorphisms,
with the third being the Beilinson-Bernstein localization theorem,
completes this proof.

\begin{align*}
\exten^{\bullet}_{U\g} (U\g\otimes_{U\n}\C,\,\,
&F_{\bullet}\Gamma(G/B,
M\otimes \lb_{\chi}))\\
&\cong\exten^{\bullet}_{U\g/\m_{\chi}}(F^*U\g\otimes_{U\n}\C,\Gamma(G/B,
M \otimes
\lb_{\chi}))\\
&\cong\exten^{\bullet}_{U\g/\m_{\chi}}(\bigoplus_{w\in W}
V(w\cdot(-2\rho-\chi)),\Gamma(G/B, M \otimes
\lb_{\chi}))\\
&\cong\bigoplus_{w\in W} \exten^{\bullet}_{D_{\chi}-mod}
(i_!\Oo_w\otimes\lb_{\chi},M\otimes \lb_{\chi})\otimes \C_{w\cdot(-2\rho-\chi)}\\
&\cong\bigoplus_{w\in W} \exten^{\bullet}_{D-mod}
(i_!\Oo_w,M)\otimes \C_{w\cdot(-2\rho-\chi)}\\
&\cong\bigoplus_{w\in W} \exten^{\bullet}_{D_{X_w}-mod}(\Oo_w,
i_w^!M)\otimes \C_{w\cdot(-2\rho-\chi)}\\
&\cong\bigoplus_{w\in W}\cohdr(X_w, i_{w}^! M)\otimes
\C_{w\cdot(-2\rho-\chi)}
\end{align*}
\end{proof}
\end{remark}

The geometric computation, in this section, of the cohomology $\cohr(\n, V)$, where $V$ is a $\g$-module that comes from a $D$-module on $G/B$, can be viewed as a recipe for reconstructing the original geometric object, namely the $D$-module, from the algebraic data of $V$ and its cohomology.  Informally, $\cohr(\n, V)$ is computed from the de Rham cohomology of the restriction of the $D$-module to the $N$-orbits.  By considering $V^g$, i.e. $g$-twists of $V$ via the adjoint action, as $g\in G$ varies, we can reconstruct the $D$-module.  In fact twisting $V$ by $g$ is equivalent, on the $D$-module side, to the pullback along the action of $g$ on $G/B$.  Thus $\cohr(\n, V^g)$ (with varying $g$) contains the data of the de Rham cohomology of the restriction of the $D$-module to the $gNg^{-1}$-orbits.  This is sufficient to recover the $D$-module; it is very natural in view of the fact that one of the orbits is a point, and varying $g$ allows the freedom of making this point, any point on $G/B$.

Let us be more precise in the following case that illustrates the general situation. Suppose that $V$ is a $\g$-module with the trivial central character.  By the Beilinson-Bernstein localization theorem we know that it comes from a $D$-module, i.e. we have $V=\Gamma(G/B,M)$ for some $D$-module $M$; let us recover it.  We have the usual short exact sequence
$$0\rightarrow\underline{\mathfrak{b}}\rightarrow \mathcal{O}_{G/B}\otimes \mathfrak{g}\rightarrow \tau_{G/B}\rightarrow
0$$ where $\mathcal{O}_{G/B}\otimes \mathfrak{g}$ is the action Lie
algebroid on $G/B$, the surjection onto the vector fields is the anchor map, and $\underline{\mathfrak{b}}$ is the kernel of the anchor map; let
$\underline{\mathfrak{n}}=[\underline{\mathfrak{b}},\underline{\mathfrak{b}}]$ and $\underline{\mathfrak{b}}/\underline{\mathfrak{n}}=\mathcal{O}_{G/B}\otimes\h$.  Then $H^{\ell(w_0)}(\underline{\mathfrak{n}}, \mathcal{L}_{-2\rho}\otimes
V)^\h$ is a $D$-module, and it follows from Theorem \ref{fincoh} that $$H^{\ell(w_0)}(\underline{\mathfrak{n}}, \mathcal{L}_{-2\rho}\otimes
V)^\h=M$$ and this is essentially the localization of $V$.  Ultimately one wishes to do the same in the affine setting, the problem  is that  there is no point orbit, however in principle it should still be possible to recover the $D$-module from its de Rham data.

\subsection{The affine setting.}\label{affine} Let us now deal with the affine
setting, namely we turn our attention to (right) $D$-modules on $\fl$.  Since we will be working with the affine Grassmannian $\gr$ and the affine flags $\fl$ extensively in what follows, we say a few words about them at this point.  Recall that $\K$ is the ring of Laurent series $\C((t))$ and we denote by $G(\K)$ the group parameterizing maps of the formal punctured disk $D^\times$ to $G$.  We have some natural subgroups: $G(\Oo)$ which parameterizes maps of the formal disk $D$ to $G$ and $I$ which is a subgroup of the latter that consists only of those maps whose center lands in $B\subset G$.  Then, roughly speaking, $\fl=G(\K)/I$ and its quotient $\gr$ is $G(\K)/G(\Oo)$.  This description is sufficient for following our geometric arguments, however for the reader interested in the foundations we point out that both can be given the structure of an ind-scheme of ind-finite type.  Furthermore, $\gr$ possesses factorization space structure and $\fl$ is a factorization module space over $\gr$.  We refer the reader to \cite{functorz, satake2} for the precise formulations.

In the affine setting
we have a choice in generalizing the finite dimensional situation.
We can consider the relationship between Iwahori orbits and Lie
algebra cohomology, or alternatively semi-infinite orbits and
semi-infinite cohomology. The latter is better suited to our
purposes and so we focus on it, briefly mentioning the former in the
remark at the end of the section.

We begin with some preliminary Lemmas establishing the shifts and twists that will appear later in the semi-infinite cohomology computations. The reader is strongly encouraged to refer to Sec. \ref{notation} when following the discussion below.  Let $w\in\Waff$, $w=\lam\wba$ (with $\lambda_w\in\Gamma$ and $\bar{w}\in W$), set
$\fraki_w=w\fraki w^{-1}$ then $\nk\cap \fraki_w$ is a semi-infinite
subspace of $\nk$. We are interested in computing the character of
the relative determinant $\dett=\dett (\nk\cap \fraki_w,\nO)$ as a
$\h$-module, as well as the relative dimension $\dimm=\dimm (\nk\cap
\fraki_w,\nO)$.  Recall that for a pair of semi-infinite subspaces
$U$ and $V$, $$\dett(U,V)=\dett(U/(U\cap V))\otimes\dett(V/(U\cap
V))^*$$ which makes sense since both $U/(U\cap V)$ and $V/(U\cap V)$
are finite dimensional, and similarly
$$\dimm(U,V)=\dimm(U/(U\cap V))-\dimm(V/(U\cap
V))$$ so that the relative dimension is an integer that need not be
non-negative.

\begin{lemma}\label{computation}
We have $\dett\cong\C_{\wba\cdot 0 +\crit(\lam)}$ and
$\dimm=-2\height\lam-\ell(\wba)$.
\end{lemma}

\begin{proof}
Observe that
\begin{align*}
\dett(\nk\cap &w\gO w^{-1},\nO)\\
&=\dett((\nk\cap \fraki_w)\oplus
(\nk\cap w\n^- w^{-1}) ,\nO)\\
&\cong\dett\otimes\dett(\nk\cap w\n^-
w^{-1})\\
&=\dett\otimes\dett(\n\cap \wba\n^-
\wba^{-1})\\
&\cong\dett\otimes \C_{-\wba\cdot 0}.
\end{align*}

While at the same time
\begin{align*}
\dett(\nk\cap &w\gO w^{-1},\nO)\\
&=\dett(\nk\cap \lam\gO \lam^{-1},\nO)\\
&=\dett(\lam\nO\lam^{-1},\nO)\\
&\cong\bigotimes_{\alpha>0}\C_{-\alpha(\lam)\alpha}\\
&=\C_{\crit(\lambda_w)}.
\end{align*}

So that $\dett\cong\C_{\wba\cdot 0 +\crit(\lam)}$.  Identically,
$\dimm=-2\height\lam-\ell(\wba)$.
\end{proof}

At this point we require an analogue of Lemma \ref{delta}. Recall that $\nkp$ and $\nop$ denote $\n(\K)\oplus t\h(\Oo)$ and $\n(\Oo)\oplus
t\h(\Oo)$ respectively. Let us
review  some basics of the semi-infinite cohomology of $\nkp$ (it is
somewhat simpler than the general case).  The complex
$\bigwedge^{\infty/2+\bullet}$ is obtained by taking the quotient of
the Clifford algebra $C:=C(\nkp\oplus\nkp^*)$ by the left ideal
generated by the standard semi-infinite isotropic subspace
 $L\subset\nkp\oplus\nkp^*$.
 Thus $\bigwedge^{\infty/2+\bullet}=C/(C\cdot L)$, where $L:=\nop\oplus(\nkp/\nop)^*$; note that $\bigwedge^{\infty/2+\bullet}$ is
a $\h$-module. Let $\vac$ denote the image of $1$ in the quotient.
The differential is written analogously with the usual cohomology
differential, namely it is
$-\frac{1}{2}f^\alpha_{\beta\gamma}:c^\beta c^\gamma c_\alpha:$,
where $f^\alpha_{\beta\gamma}$ are the structure constants of the
Lie algebra $\nkp$ with respect to the usual basis ($\{c_\alpha\}$
is the basis of $\nkp$, $\{c^\alpha\}$ the dual basis of $\nkp^*$).
Note that because $f^\alpha_{\beta\gamma}$ vanish whenever any two
of the indices coincide, the normal ordering $:c^\beta c^\gamma
c_\alpha:$ does not change the differential.  The situation modifies
readily to the case of coefficients in an $\nkp$-module $M$. Namely,
the differential in this case is given by $$c_\alpha\otimes
c^\alpha-\frac{1}{2}f^\alpha_{\beta\gamma}:c^\beta c^\gamma
c_\alpha\!:\,\,\in U(\nkp)\otimes C$$ acting on
$M\otimes\bigwedge^{\infty/2+\bullet}$.

The complex described above gives the standard semi-infinite
cohomology of $\nkp$, however one may choose a different
semi-infinite split of $\nkp$, i.e., a different isotropic subspace
of $\nkp\oplus\nkp^*$ and obtain a complex that way (as we will see
below there is not much difference, this is analogous to the
similarity between homology and cohomology in the finite case, see
the proof of Lemma \ref{shapiro}).

We will be particularly interested in the following. Let $w\in
\Waff$, recall that $\iw_w=w\iw w^{-1}$.  Then
$L_w:=(\nkp\cap\iw_w)\oplus(\nkp/(\nkp\cap\iw_w))^*$ is an isotropic
semi-infinite subspace of $\nkp\oplus\nkp^*$ (taking $w$ to be
identity in $\Waff$ we recover the standard case). We may now
consider the complex $\bigwedge_w^{\infty/2+\bullet}$ formed by
taking the quotient of the Clifford algebra $C$ by the left ideal
generated by $L_w$, let $\vac_w$ denote the image of $1$ in the
quotient. The differential is still given by the same formula (it is
an element of $C(\nkp\oplus\nkp^*)$ and so acts on any module over
the Clifford algebra), except now the normal ordering is taken with
respect to a different semi-infinite split. However, as is explained
above, for the special case of $\nkp$, this does not change the
differential.

Any choice of a non-zero element $v_w$ in the line
$\dett(\nkp\cap\iw_w,\nop)$ yields a map of $C$-modules as follows
\begin{align*}
\phi_{v_w}:\bigwedge\nolimits_w^{\infty/2+\bullet}&\longrightarrow
\bigwedge\nolimits^{\infty/2+\bullet}[-\dimm(\nkp\cap\iw_w,\nop)]\\
\vac_w&\mapsto v_w\vac
\end{align*} where the expression $v_w\vac$ is well defined since
$v_w\in\wedge^\bullet(\nkp/\nop)\otimes\wedge^\bullet \nop^*$.  This
map is an isomorphism of $C$-modules, and since the differential for
both modules is given by the same element of $C$, it is an
isomorphism of complexes.  Note that $\phi_{v_w}$ shifts the grading
as indicated and twists the $\h$-action.  Given a $\gk$-module $M$,
we denote by $\coh_w(\nkp,M)$ the cohomology of the complex
$M\otimes\bigwedge\nolimits_w^{\infty/2+\bullet}$ and omit the
subscript $w$ when it is the identity.  Observe that
$\coh_w(\nkp,M)$ is a $\h$-module.

\begin{lemma}\label{deltasemi}
Let $\delta_x$ be the right $D$-module of delta functions at $x\in
\fl$, $\chi-2\rho$ dominant regular, and $\kappa$ sufficiently
negative then
\begin{align*}
\coh(\nkp, &\Gamma(\fl,
\delta_x\otimes\lb_{\chi+\kappa}))\\
&\cong\lb_{\chi+\kappa}|_x\otimes \determinant(s_{\nkp}
x,\nop)\vac[\dimension(s_{\nkp} x,\nop)]
\end{align*}
where $s_{\nkp} x$ is the stabilizer in $\nkp$ of $x\in \fl$.
\end{lemma}

\begin{remark}
The meaning of \emph{sufficiently negative} in the statement of
Lemma \ref{deltasemi} and theorems below is as follows. We require
first of all that $\Gamma(\fl, \delta_x\otimes\lb_{\chi +\kappa})$
be irreducible as a $\gk$-module for any $x\in\fl$.  This reduces to
irreducibility of the Verma module $M_{-\chi,\kappa}$, which using
\cite{kk} can be shown to be irreducible whenever $$\kappa-\crit\leq
-(\chi-\rho,\theta)_0,$$ i.e., whenever $\kappa\leq
-(\chi,\theta)_0-1.$ We also need the exactness of the functor
$\Gamma(\fl, -\otimes\lb_{\chi +\kappa})$, however by \cite{bd} this
holds for the case when $M_{-\chi,\kappa}$ is irreducible, so the
above condition is sufficient for this as well. Another condition is
needed to ensure the degeneration of the spectral sequence that
allows us to consider only $D$-modules supported on a single orbit,
i.e., we need that $w\cdot(-\chi)\neq w'\cdot(-\chi)$ in $\h^*$
whenever $w\neq w'$ in $\Waff$.\footnote{The remark following
Theorem \ref{semicoh} can be used to show that this is also a
necessary condition for the functor $\Gamma(\fl, -\otimes\lb_{\chi
+\kappa})$ to be an embedding.} We point out that this is a
requirement on the \emph{orbit} of $-\chi$ under the dot $\Waff$
action; we ask that the action have no stabilizer. We can describe,
completely combinatorially, a sufficient condition for this to hold
i.e., for every $w\in W$, $w$ not identity, there is a root $\alpha$
of $\g$ such that
$\left(\alpha,\frac{\chi-\rho-w(\chi-\rho)}{\kappa-\crit}\right)_0\notin\mathbb{Z}$.
If we require that
$$\kappa-\crit< -2(\chi-\rho,\theta)_0,$$ then it suffices both for this and the previous two requirements.
\end{remark}

\begin{proof}
We begin exactly as we did in the proof of Lemma \ref{delta}. Reduce
to the case of $x=wI$, and note that $\Gamma(\fl,
\delta_{wI}\otimes\lb_{\chi+\kappa})$ is a simple module (by the Remark above).  Then, just as in the finite case, it
can be identified with an appropriate semi-induced module of Voronov.  Consequently (by the semi-infinite Shapiro Lemma of \cite{voronov}) we see that $$\coh_w(\nkp,\Gamma(\fl,
\delta_{wI}\otimes\lb_{\chi+\kappa}))\cong\lb_{\chi+\kappa}|_{wI}\otimes\mathbb{C}\vac_w.$$
Thus, according to the discussion preceding the present Lemma,
\begin{align*}
\coh(&\nkp,\Gamma(\fl,
\delta_{wI}\otimes\lb_{\chi+\kappa}))\\&\cong\lb_{\chi+\kappa}|_{wI}\otimes
\dett(\nkp\cap\iw_w,\nop)\vac[\dimm(\nkp\cap\iw_w,\nop)].\end{align*}
This completes the proof.

\end{proof}

We now prove the semi-infinite affine analogue of Theorem
\ref{fincoh2}.  Recall that
$$w\cdot\chi=\wba\cdot\chi-(\kappa-\kappa_c)(\lam)$$ is the affine
dot action of $\Waff$ on $\h^*$ corresponding to the level
$\kappa-\crit$.

\begin{prop}\label{semicohprep}
Let $M$ be a D-module on $\fl$, $S_w\subset \fl$ the $\Nk$ orbit
labeled by $w\in\Waff$, suppose $\chi-2\rho$ is dominant regular,
$\kappa$ sufficiently negative, then as $\h$-modules
\begin{align*}
\coh&(\nkp,\Gamma(\fl, M \otimes \lb_{\chi +\kappa}))\\
&\cong\bigoplus_{w \in \Waff} \cohdr(S_{w},i^{!}_{w}M)\otimes \C_{w
\cdot (-\chi)}[-2\height(\lam)-\ell(\wba)].
\end{align*}
\end{prop}

\begin{proof}
We begin by observing that we can reduce to the special case of
$M=i_{w*}M_0$ similarly to the finite dimensional case (we use that
$\kappa$ is sufficiently negative here). Now let $M$ be a $D$-module
on $S_w$. Below we construct an explicit map from the De Rham to the
semi-infinite cohomology.

Consider the short exact sequence of vector bundles on $S_w$
arising from the action of $\Nkp$ on its orbit:
$$\mathcal{S}tab_w\stackrel{\alpha}{\longrightarrow}
\Oo_{S_w}\otimes\nkp\stackrel{\beta}{\longrightarrow}
\mathcal{T}_{S_w}$$ and denote by $\lb_{det}$ the relative
determinant line bundle
$\dett(\mathcal{S}tab_w,\Oo_{S_w}\otimes\nop)$, so that we have a
natural map
$$\psi:\bigwedge\nolimits^\bullet\mathcal{T}_{S_w}\otimes\lb_{\dett}[\dimm]\rightarrow
\Oo_{S_w}\otimes\bigwedge\nolimits^{\infty/2+\bullet}\nkp$$ as in
the proof of Theorem \ref{fincoh2}.  This is known as the ``fermions canceling the determinantal anomaly". Similarly, $\psi$ extends to
$$
\widetilde{\psi}: i_{w\cdot}(M\otimes
\lb_{\kappa+\chi}|_{S_w}\otimes\bigwedge\nolimits^\bullet\mathcal{T}_{S_w}\otimes\lb_{\dett})[\dimm]
\rightarrow
i_{w*}M\otimes\lb_{\kappa+\chi}\otimes\bigwedge\nolimits^{\infty/2+\bullet}\nkp
$$ that is a morphism of complexes of sheaves on $\fl$.

Note that $\lb_{\dett}\otimes \lb_{\kappa+\chi}|_{S_w}$ is
canonically trivialized by the $\Nkp$ action contributing only a
twist by a $\h$ character $(\lb_{\dett}\otimes
\lb_{\kappa+\chi}|_{S_w})|_{w\Iw}\cong \C_{w\cdot (-\chi)}$.  This is
where we use Lemma \ref{computation}, the only difference is the
extra $t\ho$ term that does not affect the computation. There is
also a shift by $\dimm$ that was already noted in the above, thus on
the level of cohomology we have
\begin{align*}
\cohdr(S_{w},M)\otimes \C_{w \cdot
(-\chi)}&[-2\height(\lam)-\ell(\wba)]\\
&\longrightarrow\coh(\nkp,\Gamma(\fl, i_{w*}M \otimes \lb_{\chi
+\kappa})).
\end{align*}
The map above commutes with direct limits, so it is sufficient to
consider the case when $M$ is coherent with finite dimensional
support, so that $M=i_* M_0$, with $i:X\hookrightarrow S_w$ the
inclusion of a smooth finite dimensional $X$ that contains the
support (such an $X$ exists since $S_w$ is smooth). Then $M_0$ has a
finite resolution by finite sums of $\mathcal{D}_X$ and their direct
summands, and so we may assume that $M_0=\mathcal{D}_X$. In that
case both sides of the above can be considered $\Oo_X$-modules
(locally free), and the map becomes an $\Oo_X$ morphism.  It is thus
sufficient to check that it is an isomorphism on every fiber. This
reduces to checking the statement for $M_0=\delta_x$ with $x\in X$,
and that is the content of Lemma \ref{deltasemi}.

\end{proof}

One is actually interested in the BRST reduction, which has the
advantage of producing a vertex algebra if we begin with one.  The
following addresses that issue.

\begin{theorem}\label{semicoh}
Let $M$ be a D-module on $\fl$, $S_w$ the $\Nk$ orbit labeled by $w
\in \Waff$, $\pi_{\alpha}$ the irreducible $\hcrit$-module of
highest weight $\alpha$, suppose $\chi-2\rho$ is dominant regular,
$\kappa$ sufficiently negative, then as $\hcrit$-modules\footnote{The $\hcrit$-module $\pi_0$ is known as the Heisenberg vertex algebra, and its representation theory is the same as that of $\hcrit$.  Thus Theorem \ref{semicoh} is equivalently viewed as describing the BRST reduction as a $\pi_0$-module.}
\begin{align*}
\coh&(\nk,\Gamma(\fl, M \otimes \lb_{\chi
+\kappa}))\\
&\cong\bigoplus_{w \in \Waff} \cohdr(S_{w},i^{!}_{w}M)\otimes \pi_{w
\cdot (-\chi)}[-2\height(\lam)-\ell(\wba)].
\end{align*}
\end{theorem}

\begin{proof}
An analogue of the Kashiwara theorem states that the category of
$\hcrit$-modules with locally nilpotent $t\ho$ action is equivalent
to the category of $\h$-modules, with $t\ho$ invariants in one
direction and induction in the other giving the equivalence\footnote{This statement, without refering to it as Kashiwara theorem, is explained in \cite{thebook}.}.  Thus
it suffices to show that $t\ho$ acts locally nilpotently on
$\coh(\nk,\Gamma(\fl, M \otimes \lb_{\chi +\kappa}))$, or in light
of the above we may assume that $M=i_{w*}i_w^! M$ for a fixed
$w\in\Waff$. As before we may reduce this to the case
$M=i_{w*}i_*\mathcal{D}_X$, making sure that $w\Iw\in X$.  So that
as a $t\ho\otimes\Oo_X$-module $$\coh(\nk,\Gamma(\fl,
i_{w*}i_*\mathcal{D}_X \otimes \lb_{\chi
+\kappa}))\cong\Oo_X\otimes\coh(\nk,\Gamma(\fl, \delta_{w\Iw}
\otimes \lb_{\chi +\kappa})).$$  This is sufficient since the $t\ho$
action on both $\Gamma(\fl, \delta_{w\Iw} \otimes \lb_{\chi
+\kappa})$ and $\bigwedge^{\infty/2+\bullet}\nk$ is locally
nilpotent.

\end{proof}

\begin{remark}
Considering instead Iwahori orbits and the usual cohomology one has
the formula: $$\cohr(\iw^+,\Gamma(\fl, M\otimes\lb_{\kappa+\chi}))
\cong\bigoplus_{w\in \Waff}\cohdr(X_w,i^{!}_w M)\otimes
\C_{w\cdot(-\chi)}[-\ell(w)]$$ as $\h$-modules.  This can be shown using
averaging (thus reducing the general problem to the case of constant
$D$-modules on orbits that correspond to co-Verma modules).  A proof
of Proposition \ref{semicohprep} can then be extracted from the
consideration of the above formula for an appropriate sequence of
Iwahori conjugates. This is the approach suggested by A. Beilinson
and D. Gaitsgory and followed in \cite{shapiro}.
\end{remark}

\begin{remark}
If $M$ is a right $D$-module on $G/B$ and
$i:G/B\hookrightarrow\fl$ is the inclusion of the fiber of
$p:\fl\twoheadrightarrow\gr$ over $\gt$ , then the natural map

$$\cohr(\n,\Gamma(G/B,M\otimes\lb_{\chi}))
\longrightarrow \coh(\nk,\Gamma(\fl, i_* M \otimes \lb_{\chi
+\kappa}))$$  is an isomorphism onto the highest weights.

\end{remark}

\section{The BRST reduction.}\label{brstreduction}
Let $A$ be a $\gl$-module, then by the geometric Satake isomorphism
\cite{ginsat, satake} there is a $\gt$-equivariant $D$-module $\A$
on $\gr$ such that $\cohdr(\gr,\A)=A$ (disregarding the grading, in
fact the cohomology is rarely concentrated in degree $0$). Let us
compute the BRST reduction of $\Gamma(\gr,\A\otimes\lb_{\kappa})$.
The tools are Theorem \ref{semicoh} and the \mv theorem
\cite{satake,satake2}.

\begin{prop}\label{generalbrst}
Let $A(\lambda)$ denote the $\lambda$ weight space of $A$, then as
$\hcrit$-modules
\begin{align*}
\coh(\nk,\Gamma&(\gr,\A\otimes\lb_\kappa))\\
&\cong\bigoplus_{\dind{\lambda\in\Gamma}{w\in
W}}A(\lambda)\otimes\pi_{w\cdot 0 -(\kappa-\crit)\lambda}[-\ell(w)].
\end{align*}
\end{prop}

\begin{proof}

The notation comes from the diagram below.
$$\xymatrix{\fl \ar@{->>}[d]_{p} & S_w \ar@{_{(}->}[l]_{i_w} \ar@{->>}[d]^{\tilde{p}}\\
\gr & S_{\lam} \ar@{_{(}->}[l]^{i_{\lam}} }$$

We begin  by observing that
$\Gamma(\gr,\A\otimes\lb_{\kappa})=\Gamma(\fl,p^*\A\otimes\lb_{\kappa+2\rho})$, since the fibres of $p$ are (non-canonically) $G/B$, i.e., compact; the pullback is of right $D$-modules and so we need the factor $\lb_{2\rho}$ to make sure that $p^*\A\otimes\lb_{\kappa+2\rho}$, when restricted to the fibres of $p$, is just $\Oo_{G/B}$.

To apply Theorem \ref{semicoh} we will need $\cohdr(S_w,
i^{!}_{w}p^* \A)$, while the \mv theorem tells us that
$\cohdr(S_{\lambda},
i^{!}_{\lambda}\A)=A(\lambda)[2\height\lambda]$.  Note that since
$p$ is smooth with fiber $G/B$, and $\tilde{p}$ is smooth with fiber
$X_{\wba}$, we observe that $i^{!}_{w}p^*
\A\cong\tilde{p}^*i^{!}_{\lam}\A[-\ell(w_0)+\ell(\wba)]$.  Thus
$$\cohdr(S_w, i^{!}_{w}p^* \A)\cong\cohdr(S_{\lam},
i^{!}_{\lam}\A)[-\ell(w_0)+2\ell(\wba)].$$  This together with
re-indexing, and setting $w=ww_0$ yields the result.

\end{proof}

\begin{remark}
The method of the proof above can be used to give a version of
Theorem \ref{semicoh} for the affine Grassmannian $\gr$ with $\Nk$
orbits $S_\lambda$ indexed by $\lambda\in\Gamma$. Namely, as
$\hcrit$-modules
\begin{align*}
\coh(\nk,\Gamma&(\gr,M\otimes\lb_\kappa))\\
&\cong\bigoplus_{\dind{\lambda\in\Gamma}{w\in
W}}\cohdr(S_\lambda,i^!_\lambda M)\otimes\pi_{w\cdot 0
-(\kappa-\crit)\lambda}[-2\text{ht}\lambda-\ell(w)].
\end{align*}
\end{remark}

Let $A=\Oo_{\gl}$, and call the resulting $\gk$-module $\skg$.
Observe that $\skg$ is actually a $\gl\times\gk$-module due to the
other action of $\gl$ on $\Oo_{\gl}$. Define the following
$\hl\times\hcrit$-module
$$\V^{\bullet}=\bigoplus_{\dind{\lambda\in\Gamma}{w\in
W}}\C_{-\lambda}\otimes\pi_{w\cdot
0-(\kappa-\crit)\lambda}[-\ell(w)].$$

\begin{cor}\label{hmodofreduction}
As a $\gl\times\hcrit$-module

$$\coh(\nk,\skg)\cong\Gamma(\gl/\hl,\gl\times_{\hl}\V^{\bullet}).$$

\end{cor}

\begin{proof}
By Proposition \ref{generalbrst}
$$\coh(\nk,\skg)\cong\bigoplus_{\dind{\lambda\in\Gamma}{w\in
W}}\bigoplus_{\chi\in\Gamma^+}V^*_\chi\otimes V_\chi
(\lambda)\otimes\pi_{w\cdot 0-(\kappa-\crit)\lambda}[-\ell(w)]$$ and
$\bigoplus_{\chi\in\Gamma^+}V^*_\chi\otimes V_\chi (\lambda)$
 naturally identifies with
 $\Gamma(\gl/\hl,\gl\times_{\hl}\C_{-\lambda})$.
\end{proof}

\subsection{The chiral structure.}
If we assume that $A$ above is also a unital $\gl$-equivariant
commutative algebra, then by formal considerations we see that
$\akg=\Gamma(\gr_G,\A_G\otimes\lb_\kappa)\footnote{These are just
our $\gr$ and $\A$ from before. We will soon need to distinguish
between different Grassmanians.}$ is a vertex algebra, with a
vertex subalgebra $\vkg$ coming from the unit. See Sec. \ref{chabrief}
for more details. Note that if $A=\Gamma(X,\mathcal{B})$ where $\mathcal{B}$ is a bundle of $\gl$-equivariant commutative algebras then $\akg$ also fibers over $X$ and the fibres are $\akg_x=(\mathcal{B}_x)_\kappa(\g)$.

For our purposes, it is also useful to consider $A$ as a
$\hl$-module, and via a similar procedure we obtain another vertex
algebra $\akh=\Gamma(\gr_H,\A_H\otimes\lb_{\kappa})$.

\begin{remark}
Starting with a $\hl$-algebra $\Gamma(\hl,\Oo_{\hl})$ and
proceeding as above we get the lattice Heisenberg vertex algebra.
When $A=\Oo_{\gl}$, then $\akg=\skg$, known as the chiral Hecke algebra (Sec. \ref{chabrief}).
\end{remark}

Thus the BRST reduction of $\akg$ is not only a $\hcrit$-module, but
also a vertex algebra and below we describe the vertex algebra
structure on its subalgebra $H^{\infty/2+0}(\nk,\akg)$.  First we
need a Lemma.  Let $V$ be a finite dimensional vector space,
$\lb_\determinant$ the canonical determinant factorization line
bundle on $\gr_{GL(V)}$, and $\cl$ the constant bundle with fiber
$\bigwedge^{\bullet}_{V}$.

In what follows we briefly switch to the language of factorization algebras as the constructions involved are performed most naturally in that setting.  The languages of vertex algebras, chiral algebras and factorization algebras can be used essentially interchangeably and \cite{thebook} is an excellent dictionary.  In the proof below we use effective divisors on a curve $X$ instead of the points in $X$ as the reader may be used to.  We point out that this is basically the same thing as $X$ is one dimensional (thus effective divisors are just points with multiplicities\footnote{The dependence on multiplicities is eventually eliminated in the limit as is required.}).  However, effective divisors make sense in families and this is necessary for a proper definition of factorization structure (which is essentially a description of what happens when points collide).

\begin{lemma}\label{factorizationlemma}
$\cl$ has factorization structure and the canonical map
$$\lb_\determinant\longrightarrow\cl$$ is compatible with
factorization structures.
\end{lemma}

\begin{remark}
It was communicated to us by A. Beilinson that the Lemma is a
special case, with $G=GL_n$, of a very general situation which makes
sense for an arbitrary reductive group $G$.  Namely, consider the
vacuum integrable representation $V$ of $G(\K)^\sim$ of level
$\kappa$.  This is naturally a vertex algebra (a quotient of the
usual Kac-Moody vertex algebra).  $V$ can be realized as the dual
vector space to the space of sections of a certain positive line
bundle on $\gr$.  This line bundle admits a canonical factorization
structure, and the dual line bundle $\mathcal{L}$ embeds naturally
into $V\otimes\Oo_{\gr}$ in a way compatible with the factorization
structures.

\end{remark}
\begin{proof}
Since $\bigwedge^{\bullet}_{V}$ is a vertex algebra\footnote{See
\cite{thebook} for instance, where the structure is given by
explicit formulas.}, $\cl$ has factorization structure. From this
description of the structure one can not see directly why the
natural map above is compatible with it. There is a construction,
due to A. Beilinson, that is very similar on the level of vector
spaces to the one in \cite{thebook}, but which very naturally (i.e.,
without formulas) produces a factorization structure. Almost
tautologically this factorization algebra, call it $\Lambda$,
contains the determinant bundle as a factorization subbundle. Below
we outline the construction and show that this natural factorization
algebra is in fact the usual semi-infinite Clifford module vertex
algebra.

To define $\Lambda$ as a factorization algebra on a curve $X$, we
need to assign to each effective divisor $D$ on $X$ a vector space
$\Lambda_D$ such that when $D$ varies, $\Lambda_D$ becomes a vector
bundle (of infinite rank) on the parameter space. Furthermore, we
need to exhibit the factorization isomorphisms, i.e., for
$D=D_1+D_2$ with $D_1, D_2$ having disjoint support, we must
naturally identify $\Lambda_D$ with
$\Lambda_{D_1}\otimes\Lambda_{D_2}$.

Fix an effective divisor $D$, for $n\geq0$ let
$$W_n=V\otimes\Gamma(X,\Oo_X(nD)/\Oo_X(-nD))$$ and
$$W^*_n=V^*\otimes\Gamma(X,\omega_X(nD)/\omega_X(-nD)),$$ where
$W_n$ and $W^*_n$ are in fact non-degenerately paired via the
residue pairing.  Let $V_n=W_n\oplus W^*_n$ with its natural
bilinear form $\pair$. Note that for $m>n$, $V_n$ is naturally a
sub-quotient of $V_m$ and denote by $S_{m,n}$ the subspace of
$V_m$ that projects onto $V_n$. Let $K_{m,n}$ be the kernel of
this projection and observe that $(K_{m,n},S_{m,n})=0$.  Note that
$$A_n=S_{n,0}=V\otimes\Gamma(X,\Oo_X/\Oo_X(-nD))\oplus
V^*\otimes\Gamma(X,\omega_X/\omega_X(-nD))$$ is an isotropic
subspace of $V_n$, $A_m\subset S_{m,n}$ projects onto $A_n$, and
$K_{m,n}\subset A_m$.  Let $$\Lambda_n=C(V_n)\otimes_{\wedge
A_n}\C,$$ where $C(V_n)$ is the Clifford algebra of $V_n$. Observe
that $\Lambda_n$ is graded by assigning elements of $W_n, W_n^*$
degrees $-1$ and $1$ respectively.  Note that by above, for $m>n$,
we have $\Lambda_n\hookrightarrow\Lambda_m$ as graded vector
spaces.  Finally, $$\Lambda_D:=\varinjlim\Lambda_n$$ and one
immediately checks that it has all the properties we needed for a
factorization structure.  Namely, as $D$ varies, $V_n$, $A_n$ and
thus $\Lambda_n$ form finite dimensional vector bundles on the
parameter space.  Furthermore, a decomposition of $D$ into
disjoint $D_1$ and $D_2$ decomposes $V_n$ and $A_n$ into a direct
sum, thus $\Lambda_n$ into a tensor product.  Finally,
$\Lambda_{sD}=\Lambda_D$ for $s>0$.

The pullback of $\Lambda$ to $\gr_{GL(V)}$ naturally contains
$\lb_\determinant$ as a factorization subbundle.  Namely, for a $D$
as above, let $M\in\gr^D_{GL(V)}$, i.e., $M$ is a vector bundle on
$X$ equipped with $M|_{X-supp\, D}\cong V\otimes\Oo_X|_{X-supp\,
D}$. Thus for $n\gg 0$, $$V\otimes\Oo_X(-nD)\subset M\subset
V\otimes\Oo_X(nD)$$ and denote by $L_M$ the image of
$\Gamma(X,M/\Oo_X(-nD))$ in $W_n$.  Then $L_M\oplus L_M^\perp\subset
V_n$ is an isotropic subspace, and let $\ell_M$ be the line in
$\Lambda_n$ annihilated by $L_M$.  Then the image of $\ell_M$ in
$\Lambda_D$ is naturally identified with $\lb_\determinant|_M$.  One
immediately sees that the factorization isomorphisms are compatible.

It remains to show that $\Lambda$ is isomorphic to $\cl$. First,
observe that they are naturally isomorphic as vector spaces by
construction.  Second, choose a torus $H\subset GL(V)$ and restrict
the above factorization compatible map to $\gr_H$, i.e., we have
$$\lb_\determinant\rightarrow\Lambda\otimes\Oo_{\gr_H}.$$ Let
$\delta$ be the $D$-module of delta functions at every closed point
of $\gr_H$.  Applying $\Gamma(\gr_H,-\otimes\delta)$ to the above,
we obtain a map of factorization algebras on $X$ from a lattice
Heisenberg to $\Lambda\otimes\Oo_{\mathcal{J}\Check{H}}$, where
$\Oo_{\mathcal{J}\Check{H}}$ is the commutative factorization
algebra of functions on the jet scheme of the dual torus.  Composing
with the restriction to $1\in\mathcal{J}\Check{H}$ we obtain the
usual boson-fermion correspondence on the level of vector spaces.
Since the map is compatible with factorization structure, we are
done.
\end{proof}

\begin{remark}
For any $\kappa$, there is a natural map of factorization bundles
on $\gr_H$:
$$\lb_\kappa\rightarrow\Gamma(\gr_H,\lb_\kappa\otimes\delta)\otimes\Oo_{\gr_H}$$
obtained by taking the dual of
$\lb_\kappa^*\leftarrow\Gamma(\gr_H,\lb_\kappa^*)\otimes\Oo_{\gr_H}$.
Applying $\Gamma(\gr_H,-\otimes\delta)$ to it, we obtain the
co-action map that is the essence of the definition of the lattice
Heisenberg according to \cite{cha}.
\end{remark}

Equipped with the above we can proceed.

\begin{prop}\label{zeropartbrstprop}
As vertex algebras
$$H^{\infty/2+0}(\nk,\akg)\cong A_{\kappa-\crit}(\h).$$
\end{prop}

\begin{proof}
Consider the diagrams below. On the left $i$ is the inclusion, $p$
the usual projection and $a$ the adjoint action map, on the right
are the induced maps on the corresponding Grassmannians:
$$\xymatrix{G & B\ar[r]^a\ar@{->>}[d]^p\ar@{_{(}->}[l]_i & GL(\n)& \gr_G & \gr_B\ar[r]^a
\ar@{->>}[d]^p\ar@{_{(}->}[l]_i & \gr_{GL(\n)}\\
& H & & & \gr_H &}$$ and everything is compatible with the
factorization structure.  Call $\A_G$ the $D$-module on $\gr_G$
corresponding to $A$ under the Satake transform, denote by $\A_H$
the one on $\gr_H$.  We have the level bundle $\lb_\kappa$ on
$\gr_G$, and $\lb_\determinant$ the canonical determinant line
bundle on $\gr_{GL(\n)}$.  Then by the \mv theorem $p_*
i^!\A_G\cong\A_H$. (There are cohomology shifts appearing in the \mv
theorem, but they simply compensate for the modified commutativity
constraint. The statement should be interpreted to mean an
isomorphism of factorization sheaves.) Thus $\Gamma(\gr_H, p_*
i^!(\A_G\otimes\lb_\kappa))\cong\akh$ as vertex algebras
($\lb_\kappa$ is trivialized along the fibers of $p$).

Note that $a^*\lb_\determinant$ is simply the line bundle on $\gr_B$
of relative determinants of the stabilizers in $\nk$ of points in
$\gr_B$.  Denote also by $\cl$ the constant bundle on $\gr_B$ with
fiber $\bigwedge^{\bullet}_{\n}$.  As is mentioned above $\cl$ has
factorization structure, furthermore $a^*\lb_\determinant$ sits
inside as a factorization sub-bundle by Lemma
\ref{factorizationlemma}. We have the following geometric version of
the map of Proposition \ref{semicohprep} $$i_\cdot (DR^\bullet_p
i^!\A_G\otimes i^* \lb_\kappa\otimes
a^*\lb_\determinant)\longrightarrow i_*
i^!\A_G\otimes\lb_\kappa\otimes\cl.$$ It is also compatible with the
factorization structure.  Applying $\Gamma(\gr_G,-)$ and taking the
cohomology gives (by Proposition \ref{generalbrst}) the desired
isomorphism
$$A_{\kappa-\crit}(\h)\stackrel{\sim}{\rightarrow}H^{\infty/2+0}(\nk,\Gamma(\gr_G,i_*
i^!\A_G\otimes\lb_\kappa)).$$ We observe that the factorization
algebra $\A_G\otimes\lb_\kappa$ is filtered with the associated
graded algebra $i_* i^!\A_G\otimes\lb_\kappa$.  As before, the
$\hcrit$ action on the reduction ensures that they have the same
vertex algebra structure on their respective cohomologies thus
completing the proof.
\end{proof}

Denote by $V_{\Gamma,\kappa-\crit}$ the unique up to isomorphism
lattice Heisenberg vertex algebra associated to the lattice
$\Gamma$ and the bilinear pairing $\pair_{\kappa-\crit}$, then we
have the following description of the 0th part of the BRST
reduction.  See Sec. \ref{chabrief} for the discussion of the chiral Hecke algebra $\skg$, in particular its description as an explicit vector space.

\begin{cor}\label{hwaofzeropart}As vertex algebras
$$H^{\infty/2+0}(\nk,\skg)\cong
\Gamma(\gl/\hl,\gl\times_{\hl}V_{\Gamma,\kappa-\crit}).$$
\end{cor}

\begin{proof}
By the preceding Theorem we need to describe the vertex algebra  $A_{\kappa-\kappa_c}(\h)$ for $A=\Oo_{\gl}$.  However as a $\hl$-equivariant algebra $\Oo_{\gl}=\Gamma(\gl/\hl, \gl\times_{\hl}\Oo_{\hl})$, i.e., it fibers over $\gl/\hl$ and we note that $A_{\kappa-\kappa_c}(\h)$ for $A=\Oo_{\hl}$ is the lattice Heisenberg vertex algebra $V_{\Gamma,\kappa-\crit}$.  Thus for $A=\Oo_{\gl}$, we have that $A_{\kappa-\kappa_c}(\h)$ also fibers over $\gl/\hl$ and $A_{\kappa-\kappa_c}(\h)=\Gamma(\gl/\hl, \gl\times_{\hl} V_{\Gamma,\kappa-\crit})$.
\end{proof}

Recall that if $A$ has a unit then $\vkg\subset\akg$, and so to
describe the vertex algebra structure on the reduction of $\akg$
one must at least understand $$\Pi:=\coh(\nk,\vkg)$$ as a vertex algebra.
It follows from Proposition \ref{generalbrst} that
$$\Pi \cong\bigoplus_{w\in W}\pi_{w\cdot 0}[-\ell(w)]$$ as
a $\hcrit$-module.  This determines the vertex algebra structure
modulo the understanding of multiplication on the highest weights.
(At this point a detour through \ref{hwaappendix} is recommended.)
These are represented in the cohomology by the cocycles $\left|
w\cdot 0\right\rangle=v_k\otimes (\omega_w)_0 \vac$, where $v_k$ and
$\vac$ are the generators of $\vkg$ and $\bigwedge^{\bullet}_{\n}$
respectively, and $\omega_w$ is the cocycle in $\bigwedge^{\bullet}
\n^*$ that spans $\cohr(\n,\C)^{w\cdot 0}$, i.e.,
$\omega_w=\text{det}(\n/(\n\cap\n_w))^*$.  With this in hand one
easily computes the leading coefficient (it will occur in degree
$0$) of the OPE between two highest weight vectors, and obtains the
following.

\begin{lemma}\label{hwavkg}
The highest weight algebra of $\Pi$ is $\cohr(\n,\C)$.
\end{lemma}

We are now able to completely describe the vertex algebra structure
on $\coh(\nk,\skg)$.  Recall that
$$\V^{\bullet}=\bigoplus_{\dind{\lambda\in\Gamma}{w\in
W}}\C_{-\lambda}\otimes\pi_{w\cdot
0-(\kappa-\crit)\lambda}[-\ell(w)]$$ can be given the structure of a
vertex algebra as described in Sec. \ref{hwaappendix}.

\begin{theorem}\label{maintheorem}As a vertex algebra
$$\coh(\nk,\skg)\cong\Gamma(\gl/\hl,\gl\times_{\hl}\V^{\bullet}).$$
\end{theorem}

\begin{remark}
We note that while the above Theorem addresses the BRST reduction of the untwisted chiral Hecke algebra $\skg$, it is readily applied to the twisted case.  More precisely, recall that $\skg_\phi$ denotes the twist of $\skg$ by a $\gl$-local system $\phi$ on $\text{Spec}(\K)$.  Then $$\coh(\nk,\skg_\phi)\cong\Gamma(\gl/\hl,\gl\times_{\hl}\V^{\bullet})_\phi$$ where the right-hand side denotes the  twist of $\Gamma(\gl/\hl,\gl\times_{\hl}\V^{\bullet})$ by a $\gl$-local system $\phi$.
\end{remark}

\begin{proof}
By Corollary \ref{hwaofzeropart} and the remark in
\ref{hwaappendix}, we see that $\Gamma(\gl/\hl,\Oo_{\gl/\hl})$ is
central in $\mathcal{H}:=\coh(\nk,\skg)$.  Thus we can realize this vertex
algebra as global sections of a sheaf of vertex algebras over
$\gl/\hl$.  Recall that Corollary \ref{hwaofzeropart} identifies
$\mathcal{H}_0:=H^{\infty/2+0}(\nk,\skg)$ as a vertex algebra with
$\Gamma(\gl/\hl,\gl\times_{\hl}V_{\Gamma,\kappa-\crit})$, and
consider for every $w\in W$ the $\mathcal{H}_0$-submodule of $\mathcal{H}$ generated by
$\left| w\cdot 0\right\rangle$; denote it by $\mathcal{H}_w$.

Irreducible representations of $V_{\Gamma,\kappa-\crit}$ are
parameterized by $\Check{\Gamma}/(\kappa-\crit)(\Gamma)$ (see for
example \cite{thebook}).  More precisely, if
$\alpha\in\Check{\Gamma}$, then the irreducible representation
indexed by $\overline{\alpha}\in
\Check{\Gamma}/(\kappa-\crit)(\Gamma)$, let us call it the highest
weight representation of highest weight $\alpha$ and denote it by
$U_\alpha$, is given as a $\hcrit$-module by
$\bigoplus_{\lambda\in\Gamma}\pi_{\alpha-(\kappa-\crit)\lambda}$; it
is generated as a $V_{\Gamma,\kappa-\crit}$-module by
$\left|\alpha\right>$.  Observe that by considering $\alpha$ instead
of $\overline{\alpha}$, $U_\alpha$ is naturally a $\hl$-module,
i.e., the $\hl$-equivariant irreducible representations of
$V_{\Gamma,\kappa-\crit}$ are indexed by $\Check{\Gamma}$ itself
(more precisely, $\left|\alpha\right>$ is the highest weight of the
$\pi_0$-module $U_\alpha^{\hl}$). We note that for $\kappa$
sufficiently negative in our sense, all the $\overline{w\cdot
0}\in\Check{\Gamma}/(\kappa-\crit)(\Gamma)$ are distinct for
different $w\in W$, thus indexing non-isomorphic representations of
$V_{\Gamma,\kappa-\crit}$.

Thus, using Corollary \ref{hmodofreduction}, $\mathcal{H}_w$ can be identified
with $\Gamma(\gl/\hl,\gl\times_{\hl} U_{w\cdot 0})$, and
$\mathcal{H}=\bigoplus_{w\in W}\mathcal{H}_w$ as an $\mathcal{H}_0$-module.  Let
$A^\lambda=\bigoplus_{\chi\in\Gamma^+}V^*_\chi\otimes
V_\chi(\lambda)$ so that $\mathcal{H}=\bigoplus_{\lambda,
w}A^\lambda\otimes\pi_{w\cdot 0-(\kappa-\crit)\lambda}$.  Then
$$hwa(\mathcal{H})=\bigoplus_{\lambda, w}A^\lambda\otimes\left|\lambda+w\cdot
0\right>,$$ where we retain the $\left|\lambda+w\cdot 0\right>$ to
keep track of the a priori different $A^\lambda$. Knowledge of $\mathcal{H}$
as an $\mathcal{H}_0$-module allows us to compute (for $a_\lambda\in
A^\lambda$, $a_\chi\in A^\chi$ and $w,w'\in W$ such that
$\omega_w\cdot\omega_{w'}=\pm\omega_{w''}$):
\begin{align*}
a_\lambda\otimes&\left|\lambda+w\cdot 0\right>\cdot
a_\chi\otimes\left|\chi+w'\cdot 0\right>\\
&=a_\lambda\otimes\left|\lambda\right>\cdot 1\otimes\left|w\cdot
0\right>\cdot
a_\chi\otimes\left|\chi\right>\cdot 1\otimes\left|w'\cdot 0\right>\\
&=(-1)^{\ell(w)(\chi,\chi)+w\cdot 0(\chi)}a_\lambda\otimes\hla\cdot
a_\chi\otimes\hchi\cdot 1\otimes\left|w\cdot
0\right>\cdot 1\otimes\left|w'\cdot 0\right>\\
&=\pm(-1)^{\ell(w)(\chi,\chi)+w\cdot 0(\chi)}a_\lambda a_\chi
\otimes\left|\lambda+\chi\right>\cdot 1\otimes\left|w''\cdot
0\right>\\
&=\pm(-1)^{\ell(w)(\chi,\chi)+w\cdot 0(\chi)}a_\lambda a_\chi
\otimes\left|\lambda+\chi+w''\cdot 0\right>.
\end{align*} We conclude that $hwa(\mathcal{H})\cong hwa(\mathcal{H}_0)\widetilde{\otimes}\,\cohr(\n,\C)$
and the claim follows.

\end{proof}

As was mentioned in the introduction, the unramified case of the geometric local Langlands correspondence manifests itself in our situation in the form of the $D$-modules on the affine flags that we called monodromy annihilators.  Recall that such a $D$-module $M$ has the property that the monodromy action on $\mathcal{Z}(V)$, with $V$ any representation of $\gl$, becomes trivial on $\mathcal{Z}(V)\star M$.  The importance of this notion for us is that these $M$ provide $\gl$-equivariant representations of the untwisted chiral Hecke algebra $\skg$ via $\Gamma(\fl,(\mathcal{Z}(\Oo_{\gl})\star-)\otimes\lb_{\kappa+\chi})$.  In fact, conjecturally, these are all of them.

In particular, $D$-modules pulled back to the affine flags from the affine Grassmannian are in some sense the most important examples of the monodromy annihilators.\footnote{Their importance, conjectural and otherwise, is discussed in the introduction.}  To obtain a series of $(\skg,\gl)$-modules from them one need not even leave the affine Grassmannian.  Recall that for a $D$-module $M$ on $\gr$, we have that $\Gamma(\gr,(\widetilde{\Oo}_{\gl}\star M)\otimes\lb_\kappa)$ is an $(\skg,\gl)$-module.\footnote{This is the same $(\skg,\gl)$-module as the one obtained via pullback to $\fl$, i.e. it is isomorphic to $\Gamma(\fl,(\mathcal{Z}(\Oo_{\gl})\star \pi^* M)\otimes\lb_{\kappa+2\rho})$.  Note that in order to obtain all of the $(\skg,\gl)$-modules that come from $\gr$ one does need to pull back to $\fl$ first as otherwise any twist other than by $2\rho$ is unavailable.}  We would like to consider its BRST reduction and describe it as a module over the BRST reduction of  $\skg$ itself.

It follows from Theorem \ref{maintheorem} that the BRST reduction of a $\gl$-equivariant $\skg$-module $V$ fibers equivariantly over $\gl/\hl$ so that it is completely determined by the structure of the fiber over $1\in\gl/\hl$, let us denote it by $\mathcal{B}(V)$, as a $\hl\times\gm$-equivariant $\V^{\bullet}$-module.  This itself is determined by the structure of $\mathcal{B}(V)^{\hl}$ as a $\gm$-equivariant, i.e. graded, $\Pi$-module.  In the case under consideration $$V=\Gamma(\gr,(\widetilde{\Oo}_{\gl}\star M)\otimes\lb_\kappa)$$ and $$\mathcal{B}(V)^{\hl}=\coh(\nk,\Gamma(\gr,M\otimes\lb_\kappa)).$$  The latter can be computed as a $\pi_0$-module using the Remark that follows Proposition \ref{generalbrst}, and as in the proof of Lemma \ref{hwavkg} we see that the action of the whole $\Pi$ is ``free".  Consequently we obtain the following Corollary.

\begin{cor}\label{conj}
Let $M$ be a $D$-module on $\gr$ and $$\mathcal{H}(M):=\bigoplus_{\chi\in\Gamma}\cohdr(S_\chi,i^!_\chi M)(\chi)[-2\text{ht}\chi]$$ the associated $\hl\times\gm$-module, then as $(\V^\bullet,\hl\times\gm)$-modules $$\mathcal{B}(\Gamma(\gr,(\widetilde{\Oo}_{\gl}\star M)\otimes\lb_\kappa))\cong \V^\bullet\otimes \mathcal{H}(M)$$ where $\V^\bullet$ is viewed as a $\hl\times\gm$-equivariant module over itself and so can be twisted by the $\hl\times\gm$-module $\mathcal{H}(M)$.
\end{cor}

One may thus conjecture that the $\gl$-equivariant $\skg$-modules that arise as $V=\Gamma(\gr,(\widetilde{\Oo}_{\gl}\star M)\otimes\lb_\kappa)$ are characterized by the property that $\mathcal{B}(V^g)$ is of the form $\V^\bullet\otimes \mathcal{M}_g$ for every $g\in G(\K)$, where $\mathcal{M}_g$ is some $\hl\times\gm$-module.  By interpreting $\mathcal{M}_g$ as $\mathcal{H}(g^* M)$ for some $D$-module $M$ on $\gr$ one should be able to recover $M$ itself.

The BRST reduction of other series of $(\skg,\gl)$-modules that come from $\gr$, i.e. those arising from twisting by a character other than $2\rho$, can be similarly described through the structure of their fibres over $1\in\gl/\hl$ as $\hl\times\gm$-equivariant $\V^{\bullet}$-modules.  They are still ``free", though now modeled not on $\V^{\bullet}$ itself, but rather on a shift of it.  This is similar and in fact caused by, a similar phenomenon that occurs for lattice Heisenberg modules; they include the lattice Heisenberg itself and a finite number of its shifts.

The situation for other monodromy annihilators on $\fl$, i.e those that do not arise as pullbacks from $\gr$, is more complicated.  After applying the BRST reduction functor one can still restrict to the fiber over $1\in\gl/\hl$, however the resulting module over $\V^\bullet$ is no longer ``free" and is in general rather arbitrary.  However, if we restrict our attention only to the $V_{\Gamma,\kappa-\crit}$-module\footnote{Recall that $V_{\Gamma,\kappa-\crit}=\V^0\subset\V^\bullet$.} structure, then the situation is again very manageable.  Namely, recall that the irreducible $(V_{\Gamma,\kappa-\crit},\hl)$-module $U_\alpha$ is characterized by the property that $U_\alpha^{\hl}\cong\pi_\alpha$ as a $\pi_0$-module.  As before we have that $$\mathcal{B}(\Gamma(\fl,(\mathcal{Z}(\Oo_{\gl})\star M)\otimes\lb_{\kappa+\chi}))^{\hl}=\coh(\nk,\Gamma(\fl, M\otimes\lb_{\kappa+\chi}))$$ and the latter can be computed as a $\pi_0$-module using Theorem \ref{semicoh}.  Let us summarize as follows.

\begin{cor}\label{conj2} Let $M$ be a monodromy annihilator $D$-module on $\fl$, and set $V=\Gamma(\fl,(\mathcal{Z}(\Oo_{\gl})\star M)\otimes\lb_{\kappa+\chi})$, then as $(V_{\Gamma,\kappa-\crit}, \hl\times\gm)$-modules
$$\mathcal{B}(V)\cong\bigoplus_{w \in \Waff} U_{w
\cdot (-\chi)}\otimes  \cohdr(S_{w},i^{!}_{w}M)[-2\height(\lam)-\ell(\wba)]$$ where $w=\lam\wba$.
\end{cor}

\section{Appendix.}\label{appendix}
Here we collect some auxiliary information that we hope will make
the paper more accessible to the reader.

\subsection{The chiral Hecke algebra.}\label{chabrief}
The chiral Hecke algebra, introduced by Beilinson-Drinfeld, is
defined using the geometric version of the Satake isomorphism
\cite{ginsat, satake}, which is an equivalence (of tensor
categories) between the category of representations of the Langlands
dual group $\gl$ and the graded (by dimension of support) category
of $\gt$-equivariant $D$-modules on the affine Grassmannian (here
the tensor structure is given by convolution). The functor from
$D$-modules to $\gl$ representations is just $\cohdr(\gr,-)$. Under
this equivalence a commutative algebra structure on any $\gl$-module
produces a chiral algebra structure on the $\lk$-twisted global
sections ($\kappa$ chosen negative integral) of the corresponding
$D$-module as follows.

Let $A$ be a commutative algebra and a $\gl$-module such that the
multiplication $m:A\otimes A\rightarrow A$ is a map of
$\gl$-modules.  Let $\A$ be the corresponding (under the Satake
isomorphism) $\gt$-equivariant $D$-module on $\gr$, and
$\tilde{m}:\A
* \A\rightarrow \A$ the corresponding map of $D$-modules.  Let $X$
be a curve and consider the diagram ($\Delta$ is the embedding of
the diagonal, $j$ of the complement):
$$\xymatrix{\gr_X\ar@{->>}[d]^p\ar@{^{(}->}[r]^-{\Delta} & \gr_X^{(2)}\ar@{->>}[d]^p &
\gr_X\times\gr_X|_U\ar@{_{(}->}[l]_-{j}\ar@{->>}[d]^p\\
X\ar@{^{(}->}[r]^-{\Delta}& X\times X & U\ar@{_{(}->}[l]_-{j}}$$
One of the definitions of $\A *\A$ is as $\Delta^!
j_{!*}(\A\boxtimes\A)|_U[1]$, and so we get the diagram below:
$$\xymatrix{0\ar[r] & j_{!*}(\A \boxtimes \A)|_U \ar[r] & j_{*}(\A
\boxtimes
\A)|_U \ar[r] \ar[rd]^\alpha & \Delta_* (\A * \A) \ar[d]^{\Delta_* \tilde{m}}\ar[r] & 0\\
& & & \Delta_* \A &}$$ On $\gr_X^{(2)}$, we have
$\lb^{(2)}_\kappa$ providing the factorization structure on the
level bundle $\lb_\kappa$. When we twist the morphism $\alpha$ in
the diagram above by $\lb^{(2)}_\kappa$, we obtain
$$j_{*}(\A\otimes\lb_\kappa \boxtimes \A\otimes\lb_\kappa)|_U
\longrightarrow\Delta_* (\A\otimes\lb_\kappa).$$ By applying
$p_\cdot$  to the above, which is exact, we get a chiral bracket on
$\Gamma(\gr,\A\otimes\lb_\kappa)$. Note the use of $\A$ for both the
$D$-module on $\gr$ and also on $\gr_X$.  We denote by $p_\cdot$ the
direct image functor on the category of $\Oo$-modules, to be
contrasted with $p_*$ playing the same role for the category of
$D$-modules.

\begin{remark}
It is worthwhile to note that if instead of $p_\cdot$ above, we
apply $p_*$, necessarily to the untwisted version of the diagram,
then we again obtain a chiral bracket, on $A$ this time, which can
be constructed in a standard way from the commutative associative
product on $A$.  Namely, in the vertex algebra language
$Y(a,z)=L_a$, for $a\in A$ and $L_a$ denoting the left
multiplication operator.  Alternatively, in the chiral algebra
language, the multiplication on $A$, gives a morphism of $D$-modules
on $X$ $$A^r\otimes^! A^r\rightarrow A^r$$ where
$A^r=A\otimes_{\C}\Omega_X$ and $\otimes^!$ denotes the tensor
product of right $D$-modules obtained from the standard $\otimes$ on
left $D$-modules via the usual right-left identification.  The
chiral bracket is then constructed in the following diagram:
$$\xymatrix{0\ar[r] & A^r\boxtimes A^r\ar[r] & j_{*}j^*A^r\boxtimes A^r
\ar[r] \ar[rd] & \Delta_* (A^r\otimes^! A^r) \ar[d]\ar[r] & 0\\
& & & \Delta_* A^r &}$$

\end{remark}

Due to the nature of the commutativity constraint, the chiral
algebra we obtain is graded.  One way to describe the grading is
to say that the component arising via the Satake isomorphism from
the $\gl$-module $V_\chi$ has parity $(\chi,\chi)_{\kappa-\crit}$
mod $2$.

This procedure applied to the trivial representation yields the
Kac-Moody vertex algebra $\vkg$, while the regular representation
produces the chiral Hecke algebra $\Hk$:
$$\Hk=\bigoplus_{\chi\in\Gamma^+}
V_\chi^*\otimes\Gamma(\gr,I_\chi\otimes\lb_\kappa)$$ where
$I_\chi=i_{!*}\Omega_\chi$ the standard $\gt$-equivariant
$D$-module supported on $\gr^\chi$.  As is mentioned above, the
parity of $V_\chi^*\otimes\Gamma(\gr,I_\chi\otimes\lb_\kappa)$ is
$(\chi,\chi)_{\kappa-\crit}$ mod $2$.

\begin{remark}
If $G=H$, i.e., $G$ is a torus, then $\skg=V_{\Gamma,\kappa}$, the
lattice Heisenberg vertex algebra.  Its representation theory, in
this case, is well understood, and so the conjecture described in
the introduction is quite obviously true.
\end{remark}

\subsection{Highest weight algebras.}\label{hwaappendix}
The purpose of this section is to explain precisely how a
$\hk$-module structure on a vertex algebra essentially determines
it, the remaining information is encoded in what we call the highest
weight algebra (which is an example of a twisted commutative
algebra).  Related notions, necessary for our purposes are
discussed. The notation is borrowed from \cite{thebook}.

\begin{definition}
A twisted commutative algebra $A$ is first of all a
$\Gamma_A$-graded unital associative super-algebra, where $\Gamma_A$
is a lattice and the parity is given by a
$p:\Gamma_A\rightarrow\mathbb{Z}_2$.  We also require the additional
structure of a symmetric bilinear pairing $\pair: \Gamma_A\otimes
\Gamma_A\rightarrow \mathbb{Q}$, with  $(\lambda,\chi)\in\mathbb{Z}$
if $S_{\lambda,\chi}\neq 0$, where
$S_{\lambda,\chi}:A^{\lambda}\otimes A^{\chi}\rightarrow
A^{\lambda+\chi}$ denotes the restriction of the multiplication in
$A$ (so that $\pair$ is essentially integral). Finally $A$ must
satisfy a $\pair$-twisted commutativity constraint, i.e the
following diagram must commute (if $S_{\lambda,\chi}\neq 0$):

$$\xymatrix{a\otimes b\ar@{|->}[d]_{\sigma} & A^{\lambda}\otimes
A^{\chi}\ar[d]_{\sigma}\ar[r]^{S_{\lambda,\chi}}
& A^{\lambda+\chi}\ar@{=}[d]\\
(-1)^{p(\lambda)p(\chi)+(\lambda,\chi)}b\otimes a & A^{\chi}\otimes
A^{\lambda}\ar[r]^{S_{\chi,\lambda}}& A^{\lambda+\chi}}$$ and we
note that the commutativity constraint forces a certain
compatibility between $\pair$ and $p$, namely if
$S_{\lambda,\lambda}\neq 0$ then $p(\lambda)=(\lambda,\lambda)$ mod
$2$, else it is extra data.
\end{definition}

Consider a $\hk$-module and conformal (super) vertex algebra
$$V=\bigoplus_{\lambda \in \Gamma_V}A^\lambda\otimes\pi_{\lambda}$$ where the lattice
$\Gamma_V$ comes with a map of abelian groups to  $\h^*$ (so that we
may treat the lattice points as elements of $\h^*$). We assume that
$A^\lambda$ are finite dimensional vector spaces and the action of
$\hk$ is trivially extended to $A^\lambda\otimes\pi_{\lambda}$ from
$\pi_{\lambda}$. Recall that $\pi_\lambda$ is the Fock
representation of $\hk$, i.e., it is the module generated by the
highest weight vector $\hla$ subject to $h_n\hla=0$ if $n>0$ and
$h_0\hla=\lambda(h_0)\hla$, where $h\in\h$ and $h_n=h\otimes t^n$.

Suppose that $\pi_0 \subset V$ (which is itself a Heisenberg vertex
algebra associated to the Heisenberg Lie algebra $\hk$) is a vertex
subalgebra of $V$ (we identify $\pi_0$ with $\pi_0\cdot 1_V\subset
V$), whose action on $V$ is compatible with that of $\hk$. Let
$a_\lambda\in A^\lambda$ and denote by $V_{a_\lambda}(w)$ the field
$Y(a_\lambda\otimes\hla,w)$ associated to $a_\lambda\otimes\hla\in
A^\lambda\otimes\pi_\lambda$. Then these fields completely determine
the vertex algebra structure of $V$. But an explicit computation
(essentially present in \cite{thebook}, explicitly in
\cite{shapiro}) shows that the fields themselves are determined up
to the operations
$$S_{\lambda,\chi}:A^\lambda\otimes A^\chi\rightarrow
A^{\lambda+\chi}$$ on $A=\bigoplus_{\Gamma_V} A^\lambda$ obtained as
follows: $S_{\lambda,\chi}(a_\lambda,a_\chi)\in A^{\lambda+\chi}$ is
the leading coefficient of the series
$Y(a_\lambda\otimes\hla,w)(a_\chi\otimes\hchi)$ in
$A^{\lambda+\chi}((w))$.  Note that there is a distinguished element
$1_A\in A^0$ obtained via $1_A\otimes\vac=1_V$.

\begin{definition}
We call $A=\bigoplus_{\Gamma_V} A^\lambda$ with the operations
$S_{\lambda,\chi}$ the highest weight algebra of $V$ and denote it
$hwa(V)$.
\end{definition}

\begin{remark}
Note that the commutative algebra $A^0\otimes\vac\subset V$ is in
the center of $V$.
\end{remark}

More precisely, let $\bar{\lambda}$ denote the image in $\h$ of
$\lambda$ under $\kappa$ (we use $\kappa$ to denote the isomorphism
induced by $\pair_{\kappa}$). For $h\in\h$ let
$b^h(w)_{-}=\sum_{n<0}h_n w^{-n-1}$ and $b^h(w)_{+}=\sum_{n>0}h_n
w^{-n-1}$ then we have the following Lemma.

\begin{lemma}\label{formula}
With $V$ as above
$$V_{a_\lambda}(w)=S_{\lambda,\bullet}(a_\lambda,-)\otimes
w^{\bullet(\bar{\lambda})}e^{\int\bla_{-}}e^{\int\bla_{+}}$$ and
$hwa(V)$ is a twisted commutative algebra.  The lattice is
$\Gamma_V$, the parity is inherited from $V$, and the pairing is
given by $(\lambda,\chi)=\chi(\bar{\lambda})$.
\end{lemma}

\begin{remark}
By starting with a twisted commutative algebra
$A=\bigoplus_{\lambda\in\Gamma_A}A^\lambda$ and equipped with a
homomorphism $\psi:\Gamma_A\rightarrow\h^*$, subject to the
compatibility condition
$(\lambda,\chi)=\psi(\chi)(\overline{\psi(\lambda)})$, we can define
a vertex algebra structure on the $\hk$-module $\bigoplus
A^\lambda\otimes\pi_\lambda$ via the formula in Lemma \ref{formula}.

\end{remark}

One can describe the lattice Heisenberg vertex algebra via this
approach, namely its highest weight algebra is constructed as
follows. Consider a commutative (forgetting the grading) algebra $A$
together with a $\Gamma_A$ grading, and $p$, $\pair$ as above. Then
assuming that $p(\lambda)=(\lambda,\lambda)$ mod $2$, we can modify
the multiplication on $A$ to get $\widetilde{A}$, a $\pair$-twisted
commutative algebra.  This procedure is very similar to the one
described in \cite{cha}.  Let us begin by choosing an ordered basis
$\mathcal{B}$ of $\Gamma_A$.  For $\lambda,\chi\in\mathcal{B}$,
define
$$r(\lambda,\chi)=
\begin{cases}
p(\lambda)p(\chi)+(\lambda,\chi) & \lambda >\chi\\
0 & \text{else}
\end{cases}$$
and extend to $\Gamma_A$ by linearity.  Then if
$S_{\lambda,\chi}:A^{\lambda}\otimes A^{\chi}\rightarrow
A^{\lambda+\chi}$, let
$$\widetilde{S}_{\lambda,\chi}=(-1)^{r(\lambda,\chi)}S_{\lambda,\chi}.$$
This gives $\widetilde{A}$ the required twisted commutative algebra
structure. For the lattice Heisenberg vertex algebra we start with
the commutative algebra $\C\Gamma$.  In our case $\Gamma$ is the
co-weight lattice and the level (in our case $\kappa-\crit$) is the
pairing $\pair$. The resulting twisted commutative algebra
$\widetilde{\C\Gamma}$ is the highest weight algebra of $V_{\Gamma,
\kappa-\crit}$.

\begin{definition}
Given two twisted commutative algebras $A$ and $B$, together with a
bilinear pairing $\pair:\Gamma\otimes \Gamma\rightarrow \mathbb{Z}$
($\Gamma=\Gamma_A\oplus\Gamma_B$) extending\footnote{In the case
that is of interest to us, this extension is not the trivial one.}
those on $\Gamma_A$ and $\Gamma_B$, we can form the twisted tensor
product $A\widetilde{\otimes}B$, again a twisted commutative
algebra, by letting $$a\otimes b \cdot a'\otimes
b'=(-1)^{p(\lambda)p(\chi)+(\lambda,\chi)}a\cdot a'\otimes b\cdot
b'$$ for $b\in B^{\lambda}$ and $a'\in A^{\chi}$.
\end{definition}

The statement of Theorem \ref{maintheorem} requires three things
from this section.  First we need the twisted commutative algebra
obtained from the lattice Heisenberg vertex algebra, it is described
above.  This is a non-degenerate example in the sense that all
$S_{\lambda,\chi}$ are non-$0$.  In fact this non-degeneracy alone
implies that up to isomorphism it is a lattice Heisenberg vertex
algebra.

Our second example is $H^\bullet(\n,\C)$, a very degenerate case,
namely we take as our lattice the weight lattice (the only non-$0$
components are the lines at $w\cdot 0$ for $w\in W$). The pairing
$\pair$ is $(\kappa-\crit)^{-1}$.  Note that whenever the product of
two elements of $H^\bullet(\n,\C)$ is non-$0$, their weights are
orthogonal with respect to $(\kappa-\crit)^{-1}$, so this does not
conflict with the essential integrality of $\pair$. The parity is
given by the cohomological degree modulo $2$. We note that the
triviality of $\pair$ is necessary because $H^\bullet(\n,\C)$ is
super-commutative.  This example comes up in Lemma \ref{hwavkg}.

Finally the twisted commutative algebra that we need in Theorem
\ref{maintheorem} is formed by taking the twisted tensor product of
the two examples above.  The extension of the pairing to the direct
sum of the weight and the co-weight lattices is done through their
natural pairing.  More precisely, $(w\cdot 0,\chi)=w\cdot 0(\chi)$,
i.e., it is truly a twisted product.  We call the resulting vertex
algebra $\V^{\bullet}$, thus
$$hwa(\V^{\bullet})\cong\widetilde{\C\Gamma}\widetilde{\otimes}\,H^\bullet(\n,\C).$$

\bigskip
{\footnotesize
  \begin{tabbing}
    \hbox to \parindent{}\=\hbox to \parindent{}\=\kill
    \>\textsc{Institut des Hautes Etudes Scientifiques, Bures-sur-Yvette, France};\\
    \>\>\textit{email:} \texttt{shapiro@ihes.fr}
  \end{tabbing}}

\end{document}